\documentclass[11pt,reqno]{amsart}

\newtheorem{proposition}{Proposition}[section]
\newtheorem{lemma}[proposition]{Lemma}
\newtheorem{corollary}[proposition]{Corollary}
\newtheorem{theorem}[proposition]{Theorem}

\theoremstyle{definition}
\newtheorem{definition}[proposition]{Definition}
\newtheorem{example}[proposition]{Example}
\newtheorem{examples}[proposition]{Examples}
\newtheorem{remark}[proposition]{Remark}
\newtheorem{remarks}[proposition]{Remarks}

\newcommand{\thlabel}[1]{\label{th:#1}}
\newcommand{\thref}[1]{Theorem~\ref{th:#1}}
\newcommand{\selabel}[1]{\label{se:#1}}
\newcommand{\seref}[1]{Section~\ref{se:#1}}
\newcommand{\lelabel}[1]{\label{le:#1}}
\newcommand{\leref}[1]{Lemma~\ref{le:#1}}
\newcommand{\prlabel}[1]{\label{pr:#1}}
\newcommand{\prref}[1]{Proposition~\ref{pr:#1}}
\newcommand{\colabel}[1]{\label{co:#1}}
\newcommand{\coref}[1]{Corollary~\ref{co:#1}}
\newcommand{\relabel}[1]{\label{re:#1}}
\newcommand{\reref}[1]{Remark~\ref{re:#1}}
\newcommand{\exlabel}[1]{\label{ex:#1}}
\newcommand{\exref}[1]{Example~\ref{ex:#1}}
\newcommand{\delabel}[1]{\label{de:#1}}
\newcommand{\deref}[1]{Definition~\ref{de:#1}}
\newcommand{\eqlabel}[1]{\label{eq:#1}}
\newcommand{\equref}[1]{(\ref{eq:#1})}

\newcommand{\Aut}{{\rm Aut}\,}

\newcommand{\Id}{{\rm Id}\,}

\def\lan{\langle}
\def\ran{\rangle}
\def\ot{\otimes}

\def\NN{{\mathbb N}}

\newcommand{\Cc}{\mathcal{C}}

\def\*C{{}^*\hspace*{-1pt}{\Cc}}
\def\text#1{{\rm {\rm #1}}}

\input xy
\xyoption {all} \CompileMatrices

\usepackage{color,amssymb,graphicx,amscd}

\begin{document}

\title[Classifying coalgebra split extensions of Hopf algebras]
{Classifying coalgebra split extensions of Hopf algebras}

\author{A. L. Agore}
\address{Faculty of Engineering, Vrije Universiteit Brussel, Pleinlaan 2, B-1050 Brussels, Belgium}
\email{ana.agore@vub.ac.be and ana.agore@gmail.com}

\author{C.G. Bontea}
\address{Faculty of Mathematics and Computer Science, University of Bucharest, Str.
Academiei 14, RO-010014 Bucharest 1, Romania \quad \& \quad
Faculty of Engineering, Vrije Universiteit Brussel, Pleinlaan 2,
B-1050 Brussels, Belgium} \email{cgabrielbontea@yahoo.com}

\author{G. Militaru}
\address{Faculty of Mathematics and Computer Science, University of Bucharest, Str.
Academiei 14, RO-010014 Bucharest 1, Romania}
\email{gigel.militaru@fmi.unibuc.ro and gigel.militaru@gmail.com}

\thanks{A.L. Agore is research fellow ''aspirant'' of FWO-Vlaanderen.
This work was supported by a grant of the Romanian National
Authority for Scientific Research, CNCS-UEFISCDI, grant no.
88/05.10.2011.}



\begin{abstract}
For a given Hopf algebra $A$ we classify all Hopf algebras $E$
that are coalgebra split extensions of $A$ by $H_4$, where $H_4$
is the Sweedler's $4$-dimensional Hopf algebra. Equivalently, we
classify all crossed products of Hopf algebras $A \# \, H_4$ by
computing explicitly two classifying objects: the cohomological
'group' ${\mathcal H}^{2} (H_4, A)$ and
$\textsc{C}\textsc{r}\textsc{p} (H_4, A) :=$ the set of types of
isomorphisms of all crossed products $A \# \, H_4$. All crossed
products $A \#H_4$ are described by generators and relations and
classified: they are parameterized by the set ${\mathcal Z}
{\mathcal P} (A)$ of all central primitive elements of $A$.
Several examples are worked out in detail: in particular, over a
field of characteristic $p \geq 3$ an infinite family of
non-isomorphic Hopf algebras of dimension $4p$ is constructed. The
groups of automorphisms of these Hopf algebras are also described.
\end{abstract}
\maketitle

\section*{Introduction}
Let $A$ and $H$ be two given groups. The extension problem of
H\"{o}lder asks for the classification of extensions of $A$ by
$H$, i.e. of all groups $E$ that fit into an exact sequence
\begin{eqnarray} \eqlabel{primulsir}
\xymatrix{ 1 \ar[r] & A \ar[r]^{i} & {E} \ar[r]^{\pi} & H \ar[r] &
1 }
\end{eqnarray}
The classical approach (\cite{AB}, \cite{R}) proves that any
extension of $A$ by $H$ is equivalent to a crossed product
extension and, if $A$ is an abelian group, the Schreier's theorem
shows that all extensions are classified by the second cohomology
group ${\mathcal H}^2 (H, A)$ \cite[Theorem 7.34]{R}. The result
remains valid in the non-abelian case: this time ${\mathcal H}^2
(H, A)$ is not a group anymore but only a pointed set
\cite[Exercise 8, pg. 86]{AB}. The first generalization of
Schreier's theorem from groups to Hopf algebras was given by
Sweedler \cite{Sw68}: if $H$ is a cocommutative Hopf algebra and
$A$ a commutative algebra the cohomology ${\mathcal H}^i (H, A)$
was introduced such that ${\mathcal H}^2(H, A)$ classifies all
cleft extensions of $A$ by $H$ \cite[Theorem 8.6]{Sw68}. The
graded case was studied in \cite{Sin}. The first obstacle in the
way of generalizing the extension problem from groups to the level
of Hopf algebras was overcome at the beginning of the 90's by
defining the notion of exact sequence of Hopf algebras. Nowadays
the unanimously accepted definition for this concept is the one
given in \cite[Definition 1.2.0]{AD} (see also \cite[Definition
3.1]{Ho} and \cite[Definition 1.5]{Sch93}). This is the context in
which the extensions of Hopf algebras were studied in a series of
papers \cite{AndCan}, \cite{AD}, \cite{By1}, \cite{By2},
\cite{MS}, \cite{Ma1}, \cite{Ma2}, etc. The tool for studying the
extension problem for Hopf algebras is the so-called \emph{cocycle
bicrossproduct} $A \, ^{\tau}{}\#_{\sigma} H$ introduced in
\cite[Theorem 2.9]{MS} and independently in \cite[Theorem
2.20]{AD}. $A \, ^{\tau}{}\#_{\sigma} H$ is the vector space $A
\ot H$ with a crossed product algebra structure and a crossed
coproduct coalgebra structure. The datum that constructs the
cocycle bicrossproduct $A \, ^{\tau}{}\#_{\sigma} H$ must satisfy
several compatibility conditions, some of them very technical
(\cite[Theorem 2.20]{AD}). \cite[Proposition 3.12]{AndCan} shows
that any cleft extension of Hopf algebras is equivalent to a
cocycle bicrossproduct extension. The classification of all cleft
extension was given in \cite[Section 3]{AD} where it is shown that
all cleft extensions of $A$ by $H$ are classified by a certain
cohomological object denoted by $H^{1}_{*} (H, A)$ \cite[Theorem
3.2.14]{AD}. This is probably the most general version of the
classical Schreier theorem known for Hopf algebras. Unfortunately,
its importance is rather a theoretical one: the explicit
description and classification of all (cleft) extensions of $A$ by
$H$ -- or equivalently of all cocycle bicrossproducts $A
^{\tau}{}\#_{\sigma} H$ - is a very difficult task for two reasons
(see \cite[Section 5.2]{AndCan} for details). On the one hand the
large number of compatibility conditions that need to be fulfilled
for constructing all cocycle bicrossproducts $A
^{\tau}{}\#_{\sigma} H$ makes the problem very difficult for a
computational approach. On the other hand, there is no efficient
cohomology theory for arbitrary Hopf algebras, similar to the one
from group theory, to make a direct description of $H^{1}_{*} (H,
A)$ possible. One of the few examples known is \cite[Lema 2.8]{GV}
where it is proved that any extension of $H_4$ by $H_4$ is
equivalent to the trivial extension $H_4 \hookrightarrow H_4 \ot
H_4 \twoheadrightarrow H_4$.

For this reason, in the present paper we deal with a special case
of Hopf algebra extensions, namely the \emph{coalgebra split
extensions}. Let $A$ and $H$ be two given Hopf algebras. A
coalgebra split extension of $A$ by $H$ is a pair $(E, \pi)$,
where $E$ is a Hopf algebra that fits into a sequence $A
\hookrightarrow E \stackrel{\pi} \to H$ such that the Hopf algebra
map $\pi: E \to H$ splits in the category of coalgebras and $A
\simeq E^{{\rm co}(H)}$. Several other types of split extensions
of Hopf algebras are studied in \cite{AMB10}, \cite{AMSt},
\cite{Sch2}. The coalgebra split extensions cover the extension
problem from the theory of groups (\exref{grLiecoalg}). Exactly as
in the group case, any coalgebra split extension of $A$ by $H$ is
equivalent to a crossed product extension $(A \#H, \pi_H)$
(\prref{descriereclef}). Thus, the classification of all coalgebra
split extensions of $A$ by $H$ is equivalent to the classification
of all crossed products $A \#_{f}^{\triangleright} \, H$
associated to all possible crossed systems of Hopf algebras $(A,
H, \triangleright, f)$. The classification will be given in two
ways: from the view point of the extension theory (that is, up to
an isomorphism of Hopf algebras that stabilizes $A$ and $H$) they
will be classified after we explicitly compute the cohomological
'group' ${\mathcal H}^{2} (H, A)$ which is the counterpart for
Hopf algebras of the second cohomology group from group theory.
The second and more general way of classifying such extensions
will be given by computing explicitly the second classifying
object: $\textsc{C}\textsc{r}\textsc{p} (H, A) :=$ the set of
types of isomorphisms of Hopf algebras of all crossed products $A
\# \, H$. There exists a canonical surjection ${\mathcal H}^{2}
(H, A) \twoheadrightarrow \textsc{C}\textsc{r}\textsc{p} (H, A)$.

The paper is organized as follows: in \seref{preli} we recall the
basic concepts related to crossed products of Hopf algebras.
\seref{sec2} contains some technical results: \thref{toatemorf}
describes the set of all morphisms between two arbitrary crossed
products of Hopf algebras which is our tool in the classification
problem as well as for computing the automorphisms group of a
given crossed product of Hopf algebras. \seref{exclasros} provides
an example of classification. More precisely, for a given Hopf
algebra $A$, the crossed systems $(A, H_{4}, \triangleright, f)$
are completely described in \thref{ex1thgen}: they are
parameterized by the set $\mathcal{Z}\mathcal{P}(A)$ of all
central primitive elements of $A$. For a large class of Hopf
algebras $A$, including the enveloping algebras of Lie algebras,
\thref{classinteza1} classifies this new family of Hopf algebras
by computing ${\mathcal H}^{2} (H_4, A)$ and
$\textsc{C}\textsc{r}\textsc{p} (H_4, A)$. The group $\Aut_{\rm
Hopf}(A \# H_4)$ is explicitly described. In \seref{exdetaliate}
we construct some explicit examples: we shall classify all crossed
products of the form $A \# \,  H_{4}$, for some specific Hopf
algebras $A$, namely for the polynomial Hopf algebra $k[Y]$ and
for two of its quotients in the case when ${\rm char} (k) = p \geq
3$. Let $k = \mathbb{F}_p (X_{1}, X_{2}, \cdots)$ be the field of
rational functions in indeterminates $\{X_{i}\}_{i \geq 1}$ over
the finite field $\mathbb{F}_{p}$. \coref{nrtip1} proves that  $k
\lan y ~|~ y^{p} = y \ran \, \# \, H_4 $ contains an infinite
family of non-isomorphic $4p$-dimensional Hopf algebras. In
particular, we construct an infinite number of types of Hopf
algebras of dimension $12$ over a field of characteristic $3$.

\section{Preliminaries} \selabel{preli}

Unless specified otherwise, all algebras, coalgebras, Hopf
algebras or tensor products are over an arbitrary field $k$. For a
coalgebra $C$, we use Sweedler's $\Sigma$-notation: $\Delta(c) =
c_{(1)}\ot c_{(2)}$, $(I\ot\Delta)\Delta(c) = c_{(1)}\ot
c_{(2)}\ot c_{(3)}$, etc. (summation understood). For a $k$-linear
map $f: H \ot H \to A$ we denote $f(g, \, h) = f (g\ot h)$, for
all $g$, $h\in H$. For all unexplained notations we refer to
\cite{M}. Let $A$ and $H$ be two Hopf algebras. A linear map $r: H
\rightarrow A$ is called \emph{cocentral} if $r$ is a morphism of
coalgebras and the following compatibility condition holds for any
$h\in H$:
\begin{equation}\eqlabel{0aacr}
r(h_{(1)}) \ot h_{(2)} = r(h_{(2)}) \ot h_{(1)}
\end{equation}
The set $CoZ(H, A)$ of all cocentral maps is a group with respect
to the convolution product \cite[pg. 338]{am1}. We denote by
$CoZ^{1}(H, A)$ the subgroup of $CoZ(H, A)$ of all cocentral maps
$r : H \to A$ such that $r(1) = 1$. A $k$-linear map
$\triangleright: H \otimes A \rightarrow A$ is called a
\textit{weak action} of $H$ on $A$ if for any $a$, $b\in A$, $h\in
H$:
\begin{eqnarray}
\eqlabel{1a}  h \triangleright 1_{A} &=& \varepsilon_{H}(h)1_{A}, \qquad  1_H \triangleright a = a  \\
\eqlabel{2}   h \triangleright (ab) &=& (h_{(1)} \triangleright a)
(h_{(2)} \triangleright b)
\end{eqnarray}
For a $k$-linear map $f: H \otimes H \rightarrow A$ and a weak
action $\triangleright: H \otimes A \rightarrow A$ we shall denote
by $A \#_{f}^{\triangleright} \, H$ the $k$-vector space $A\ot H$
with the multiplication given by
\begin{equation}\eqlabel{001}
(a \# h) \cdot (c \# g):= a (h_{(1)}\triangleright c)
f\bigl(h_{(2)} , g_{(1)}\bigl) \, \# \, h_{(3)}g_{(2)}
\end{equation}
for all $a$, $c\in A$, $h$, $g\in H$, where we denoted $a\ot h$ by
$a\# h$. The object $A \#_{f}^{\triangleright} \, H$ is called a
\textit{crossed product of Hopf algebras} if it is a Hopf algebra
with the multiplication \equref{001}, the unit $1_A \# 1_H$ and
the coalgebra structure given by the tensor product of coalgebras.
In this case $(A, H, \triangleright, f)$ is called a
\textit{crossed system of Hopf algebras} \cite[Definition
1.1]{agorecia}. If $f$ is the trivial cocycle, that is $f (g, h) =
\varepsilon (g) \varepsilon (h) 1_A$, for all $g$, $h\in H$ then
the associated crossed product $A \#_{f}^{\triangleright} \,H =
A\# H$ is the semi-direct (smash) product of Hopf algebras
\cite{Mo}. The following gives necessary and sufficient conditions
for $A \#_{f}^{\triangleright} \, H$ to be a crossed product of
Hopf algebras.

\begin{proposition}\prlabel{crossedrevised}
Let $A$, $H$ be Hopf algebras, $\triangleright: H \otimes A
\rightarrow A$ a weak action and $f: H \otimes H \rightarrow A$ a
$k$-linear map. The following are equivalent:

$(1)$ $A \#_{f}^{\triangleright} \, H$ is crossed product of Hopf
algebras;

$(2)$ $\triangleright: H \otimes A \rightarrow A$ and $f: H
\otimes H \rightarrow A$ are morphisms of coalgebras satisfying
the following compatibilities for any $a\in A$, $g$, $h$, $l\in
H$:
\begin{eqnarray}
\eqlabel{3}   f(h, 1_{H}) &=& f(1_{H}, h) =
\varepsilon_{H}(h)1_{A} \\
\eqlabel{4}   [g_{(1)} \triangleright (h_{(1)} \triangleright a)]
f (g_{(2)}, \,
h_{(2)}) &=& f(g_{(1)}, \, h_{(1)}) \bigl ( (g_{(2)} h_{(2)}) \triangleright a \bigl)  \\
\eqlabel{5}    \bigl(g_{(1)} \triangleright f(h_{(1)}, \,
l_{(1)})\bigl) f\bigl(g_{(2)}, \, h_{(2)} l_{(2)} \bigl) &=&
f(g_{(1)}, \, h_{(1)}) f(g_{(2)}h_{(2)}, \, l) \\
\eqlabel{6}    g_{(1)} \otimes g_{(2)} \triangleright a &=&
g_{(2)} \otimes g_{(1)} \triangleright a \\
\eqlabel{7}   g_{(1)} h_{(1)} \otimes f(g_{(2)}, \, h_{(2)}) &=&
g_{(2)} h_{(2)} \otimes f(g_{(1)}, \, h_{(1)})
\end{eqnarray}
In this case the antipode of $A \#_{f}^{\triangleright} \, H$ is
given by
\begin{equation}\eqlabel{antipod}
S(a \# g) := \Bigl(S_{A}\bigl[f\bigl(S_{H}(g_{(2)}), \,
g_{(3)}\bigl)\bigl] \# S_{H}(g_{(1)})\Bigl) \cdot \bigl(S_{A}(a)
\# 1_{H}\bigl)
\end{equation}
for all $a\in A$ and $g\in H$.
\end{proposition}

\begin{proof}
$(2) \Rightarrow (1)$ is \cite[Lemma 1.2.10]{AndN}, where the
crossed product is viewed as a special case of the cocycle
bicrossproduct \cite[Theorem 2.20]{AD} if we let the cocycle
cross-coproduct be the trivial one. It can be also obtained as a
special case of the unified product of \cite[Theorem 2.4, Examples
2.5(2)]{am1}.

$(1) \Rightarrow (2)$ Assume that $A \#_{f}^{\triangleright} \, H$
is a Hopf algebra with the above structures. In particular, it is
an associative algebra and hence \cite[Proposition 6.1.10]{dnr}
the compatibility conditions \equref{3}, \equref{4} and \equref{5}
holds. It remains to prove that $f$ and $\triangleright$ are
morphisms of coalgebras satisfying \equref{6} and \equref{7}.
Indeed, it follows from $\varepsilon \bigl( (1_A \# h) \cdot (1_A
\# g) \bigl ) = \varepsilon (1_A \# h) \varepsilon (1_A \# g)$
that $\varepsilon_A (f (h, g) ) = \varepsilon_H (h)
\varepsilon_H(g)$, for all $h$, $g\in H$. On the other hand, the
relation $\varepsilon \bigl( (1_A \# h) \cdot (a \# 1_H) \bigl ) =
\varepsilon (1_A \# h) \varepsilon (a \# 1_H)$ gives
$\varepsilon_A (h \triangleright a) = \varepsilon_H(h)
\varepsilon_A (a)$. Now, applying $I \ot \varepsilon_H \ot I \ot
\varepsilon_H$ to the relation $\Delta \bigl( (1_A \# g) \cdot
(1_A \# h) \bigl ) = \Delta  (1_A \# g) \Delta  (1_A \# h)$ we
obtain that $f$ is a coalgebra map and applying $\varepsilon_A \ot
I \ot I \ot \varepsilon_H$ to the same relation we obtain
\equref{7}. Similarly, applying $I \ot \varepsilon_H \ot I \ot
\varepsilon_H$ to $\Delta \bigl( (1_A \# h) \cdot (a \# 1_H) \bigl
) = \Delta (1_A \# h) \Delta (a \# 1_H)$ we obtain that
$\triangleright$ is a coalgebra map and applying $\varepsilon_A
\ot I \ot I \ot \varepsilon_H$ to the same relation we obtain
\equref{6}.
\end{proof}

Let $(A, H, \triangleright, f)$ be a crossed system of Hopf
algebras. There exist Hopf algebras morphisms
\begin{equation}\eqlabel{morcanfol}
i_A : A \to A \#^{\triangleright}_{f} \, H, \quad i_A (a) = a  \#
1_H, \qquad \pi_H : A \#^{\triangleright}_{f} \, H \to H, \quad
\pi_H (a \# h) = \varepsilon_A (a) h
\end{equation}
The crossed product $A \#^{\triangleright}_{f} \, H$ will be
viewed as a left $A$-module (resp. right $H$-comodule) via the
restriction of scalar through $i_A$ (resp. $\pi_H$), i.e. $a'
\cdot (a \# h) := a' a\#h$ and $a\#h \mapsto a \# h_{(1)} \ot
h_{(2)}$, for all $a$, $a'\in A$ and $h\in H$. The crossed product
of Hopf algebras is the construction responsible for the
description of the following type of extensions of Hopf algebras:

\begin{definition}\delabel{coalspex}
Let $A$ and $H$ be two given Hopf algebras. A \emph{coalgebra
split extension of $A$ by $H$} is a pair $(E, \pi)$ consisting of
a Hopf algebra $E$ that fits into a sequence
\begin{equation} \eqlabel{sirulcross}
A \hookrightarrow E \stackrel{\pi}{\to} H
\end{equation}
such that $\pi : E \to H $ is a morphism of Hopf algebras which
has a section as a coalgebra map and $A \simeq E ^{{\rm co}(H)} :=
\{x \in E \, | \, x_{(1)} \ot \pi (x_{(2)}) = x \ot 1 \} $.

Two coalgebra split extensions $(E, \pi)$, $(E', \pi')$ of $A$ by
$H$ are called \emph{equivalent} if there exists an isomorphism of
Hopf algebras $\psi : E \to E'$ that stabilizes $A$ and
co-stabilizes $H$, i.e. the following diagram commutes:
$$
\xymatrix {& A \ar[r]^{i} \ar[d]_{Id_{A}} & {E}
\ar[r]^{\pi}\ar[d]^{\psi} & H\ar[d]^{Id_{H}}\\
& A\ar[r]^{i'} & {E'}\ar[r]^{\pi'} & H}
$$
\end{definition}
Any crossed product $(A \#_{f}^{\triangleright} \, H, \, \pi_H)$
is a coalgebra split extension of $A$ by $H$ via the canonical
morphisms given by \equref{morcanfol}. Indeed, $i_H : H \to A
\#^{\triangleright}_{f} \, H$, $i_H (h) = 1_A \# h$ is a coalgebra
map and a section of $\pi_H$. Now, $A \#^{\triangleright}_{f} \,
H$ is a right $H$-comodule algebra via $\pi_H$, i.e. the right
$H$-coaction is given by $a \# h \mapsto a \# h_{(1)} \ot
h_{(2)}$. Then we can easily prove that
$$
(A \#^{\triangleright}_{f} \, H)^{{\rm co}(H)} = \{x \in A
\#^{\triangleright}_{f} \, H \, | \, x_{(1)} \ot \pi_H (x_{(2)}) =
x \ot 1 \} = A \# 1_H \cong A
$$

Conversely, using the theory of cleft extensions \cite{DT}, we
have the following result which reduces the classification of all
coalgebra split extensions to the one of crossed products:

\begin{proposition} \prlabel{descriereclef}
Let $A$ and $H$ be two Hopf algebras. Then any coalgebra split
extension $(E, \pi)$ of $A$ by $H$ is equivalent to a crossed
product extension $(A \#_{f}^{\triangleright} \, H, \, \pi_H)$ of
$A$ by $H$.
\end{proposition}

\begin{proof} Let $(E, \pi)$ be a coalgebra split extension of $A$
by $H$ and $\varphi : H \to E$ be a coalgebra map and a section
for $\pi$. Without loss of generality, we can assume that $\varphi
(1) = 1$ (otherwise we can replace $\varphi$ by $\varphi (1)^{-1}
\varphi$). Then $\varphi$ is invertible in convolution with the
inverse $ \varphi^{-1} = S_E \circ \varphi$ and moreover, using
that $\varphi: H \to E$ splits $\pi$, it is also a right
$H$-colinear map. Indeed, for any $h\in H$
$$
\rho_E (\varphi (h) ) = \varphi (h)_{(1)} \ot \pi ( \varphi
(h)_{(2)}) = \varphi (h_{(1)}) \ot \pi (\varphi (h_{(2)})) =
\varphi (h_{(1)}) \ot h_{(2)} = (\varphi \ot {\rm Id}) \circ
\Delta (h)
$$
Thus the right $H$-extension $E/A$ is cleft. It follows from
\cite[Theorem 11]{DT} that there exist well defined maps
\begin{eqnarray*}
\triangleright = \triangleright_{\varphi} : H \ot A \to A, \qquad
h \triangleright a &:=&
\varphi (h_{(1)}) \, a \, \varphi^{-1}(h_{(2)}) \\
\,\, f = f_{\varphi}: H\ot H \to A, \qquad f (g, \, h) &:=&
\varphi (g_{(1)}) \, \varphi (h_{(1)}) \, \varphi^{-1} (g_{(2)}
h_{(2)})
\end{eqnarray*}
such that $\triangleright_{\varphi}$ is an weak action and
$$
\psi : A \#^{\triangleright}_{f} \, H \to E, \qquad \psi (a \# h)
:= a \, \varphi (h)
$$
is an isomorphism of associative unitary algebras. But there is
more: $\psi : A \ot H \to E$ is also a morphism of coalgebras, as
a composition of such maps. Thus, $\psi : A
\#^{\triangleright}_{f} \, H \to E$ is an isomorphism of algebras
and coalgebras between $A \#^{\triangleright}_{f} \, H$ and $E$.
Since $E$ is a Hopf algebra we obtain that $A
\#^{\triangleright}_{f} \, H$ is in fact a Hopf algebra and $\psi
: A \#^{\triangleright}_{f} \, H \to E$ is an isomorphism of Hopf
algebras. The fact that $\psi$ stabilizes $A$ and co-stabilizes
$H$ is straightforward.
\end{proof}

\begin{example} \exlabel{grLiecoalg}
Let $A$ and $H$ be two groups. Then a Hopf algebra $E$ is a
coalgebra split extension of $k[A]$ by $k[H]$ if and only if $E
\cong k[G]$, for a group $G$ which is an extension of $A$ by $H$.
Indeed, using \prref{descriereclef}, we obtain that $E$ is a
coalgebra split extension of $k[A]$ by $k[H]$ if and only if $E
\cong k[A] \# k[H]$, for some crossed system of Hopf algebras
$(k[A], k[H], \triangleright, f)$. Now, such crossed systems of
Hopf algebras are in bijection to the usual crossed system of
groups and the bijection is given such that there exists a
canonical isomorphism of Hopf algebras $k[A] \# k[H] \cong k [A \#
H]$ (see \cite[Examples 1.2]{agorecia} for details), where $A \#H$
is a crossed product of groups, i.e. an extension of $A$ by $H$.
\end{example}

\section{Morphisms between crossed products} \selabel{sec2}
First of all we shall prove a technical result that will be our
tool in the classification of all crossed products as well as for
computing the automorphisms group of a given crossed product of
Hopf algebras.

\begin{theorem}\thlabel{toatemorf}
Let $(A, H, \triangleright, f)$ and $(A', H', \triangleright',
f')$ be two crossed systems of Hopf algebras. Then there exists a
bijective correspondence between the set of all morphisms of Hopf
algebras $\psi : A \#^{\triangleright}_{f} \, H \to A'
\#^{\triangleright'}_{f'} \, H' $ and the set of all quadruples
$(u, p, r, v)$, where $p: A \to H'$ is a morphism of Hopf
algebras, $u: A \to A'$, $r: H \rightarrow A'$ and $v: H
\rightarrow H'$ are unitary morphisms of coalgebras satisfying the
following compatibility conditions:
\begin{enumerate}
\item[(CP1)] $u(a_{(1)}) \ot p(a_{(2)}) = u(a_{(2)}) \ot p(a_{(1)})$\\
\item[(CP2)] $r(h_{(1)}) \ot v(h_{(2)}) = r(h_{(2)}) \ot v(h_{(1)})$\\
\item[(CP3)] $u(ab) = u(a_{(1)}) \, \bigl( p (a_{(2)})
\triangleright' u(b_{(1)})
\bigl) f' \, \bigl(p(a_{(3)}), \, p(b_{(2)})\bigl)$\\
\item[(CP4)] $v(h) \, v(g) = p\bigl(f(h_{(1)}, \, g_{(1)})\bigl) \, v(h_{(2)}g_{(2)})$\\
\item[(CP5)] $v(h) \, p(a) = p(h_{(1)} \triangleright a) \,
v(h_{(2)})$\\
\item[(CP6)] $r(h_{(1)}) \bigl(v(h_{(2)}) \triangleright'
r(g_{(1)})\bigl)
\, f' \, \bigl(v(h_{(3)}), \, v(g_{(2)})\bigl) = \\
u \bigl(f(h_{(1)}, \, g_{(1)})\bigl) \Bigl(p\bigl(f(h_{(2)}, \,
g_{(2)}\bigl) \triangleright' r(h_{(4)} g_{(4)})\Bigl) f' \,
\Bigl(p\bigl(f(h_{(3)}, g_{(3)}),\,
v(h_{(5)}g_{(5)})\bigl)\Bigl)$ \\
\item[(CP7)] $r(h_{(1)}) \bigl(v(h_{(2)}) \triangleright'
u(a_{(1)})\bigl)
\, f' \, \bigl(v(h_{(3)}), \, p(a_{(2)})\bigl) = \\
u(h_{(1)} \triangleright a_{(1)}) \, \bigl(p(h_{(2)}
\triangleright a_{(2)}) \triangleright' r(h_{(4)})\bigl) f' \,
\bigl(p(h_{(3)} \triangleright a_{(3)}), \, v(h_{(5)})\bigl) $
\end{enumerate}
for all $a$, $b \in A$, $g$, $h \in H$. Under the above bijection
the morphism of Hopf algebras $\psi: A \#^{\triangleright}_{f} H
\to A' \#^{\triangleright'}_{f'}  H' $ corresponding to $(u, p, r,
v)$ is given by:
\begin{equation}\eqlabel{morfcros}
\psi(a \# h) = u(a_{(1)}) \, \bigl( p(a_{(2)}) \triangleright'
r(h_{(1)}) \bigl)f' \, \bigl(p(a_{(3)}), \, v(h_{(2)})\bigl) \,\,
\#' \, p(a_{(4)}) \, v(h_{(3)})
\end{equation}
for all $a \in A$ and $h\in H$.
\end{theorem}

\begin{proof}
Let $\psi: A \#^{\triangleright}_{f} H \to A'
\#^{\triangleright'}_{f'} H' $ be a morphism of Hopf algebras. We
define
$$
\alpha : A \to A' \#^{\triangleright'}_{f'} H', \quad \alpha (a)
:= \psi (a \# 1_H), \quad \beta : H \to A'
\#^{\triangleright'}_{f'} H', \quad \beta (h) := \psi (1_A \# h)
$$
Then $\alpha : A \to A' \ot  H'$ and $\beta : H \to A' \ot H'$ are
unitary morphisms of coalgebras as compositions of such maps and
\begin{equation}\eqlabel{ecu3}
\psi (a \# h) = \psi ( (a \# 1_H) (1_A \# h) ) = \psi (a \# 1_H)
\psi (1_A \# h) = \alpha (a) \, \beta(h)
\end{equation}
for all $a\in A$ and $h\in H$. It follows from \cite[Lemma
2.1]{ABM1} that there exist four coalgebra maps $u: A \to A'$, $p:
A \to H'$, $r: H \rightarrow A'$, $v: H \rightarrow H'$ such that
\begin{equation}\eqlabel{ecu33}
\alpha (a) =  u (a_{(1)}) \ot p (a_{(2)}), \qquad \beta (h) = r
(h_{(1)}) \ot v (h_{(2)})
\end{equation}
and the pairs $(u, p)$ and $(r, v)$ satisfy the symmetry
conditions (CP1) and (CP2). Explicitly $u$, $p$, $r$ and $v$ are
defined by
$$
u(a) = (({\rm Id}\ot \varepsilon_{H'}) \circ \psi ) (a \# 1_H),
\quad p(a) = ((\varepsilon_{A'} \ot {\rm Id}) \circ \psi ) (a \#
1_H)
$$
$$
r(h) =  (({\rm Id}\ot \varepsilon_{H'}) \circ \psi ) (1_A \# h),
\quad v(h) = ((\varepsilon_{A'} \ot {\rm Id}) \circ \psi ) (1_A \#
h)
$$
for all $a\in A$ and $h\in H$. All these maps are unitary
coalgebra maps. Now, for any $a\in A$ and $h\in H$ we have:
\begin{eqnarray*}
\psi (a \# h) &=& \alpha (a) \beta(h) = \bigl( u (a_{(1)}) \#' p
(a_{(2)}) \bigl) \cdot \bigl( r
(h_{(1)}) \#' v (h_{(2)}) \bigl)\\
&=& u(a_{(1)}) \, \bigl( p(a_{(2)}) \triangleright' r(h_{(1)})
\bigl)f' \, \bigl(p(a_{(3)}), \, v(h_{(2)})\bigl) \,\, \#' \,
p(a_{(4)}) \, v(h_{(3)})
\end{eqnarray*}
i.e. \equref{morfcros} also holds. Thus any bialgebra map $\psi: A
\#^{\triangleright}_{f} H \to A' \#^{\triangleright'}_{f'} H' $ is
determined by the formula \equref{morfcros}, for some unique
quadruple of unitary coalgebra maps $(u, p, r, v)$.

Now, we prove that a map $\psi$ given by \equref{ecu3} is a
morphism of algebras if and only if $\alpha : A \to A'
\#^{\triangleright'}_{f'} H'$ is an algebra map and the following
compatibility conditions hold:
\begin{equation}\eqlabel{ecu111}
\beta(h) \, \beta(g) = \alpha\bigl(f(h_{(1)}, \, g_{(1)})\bigl)\,
\beta(h_{(2)}g_{(2)})
\end{equation}
\begin{equation}\eqlabel{ecu11}
\beta (h) \, \alpha (b) = \alpha (h_{(1)} \triangleright b )  \,
\beta (h_{(2)})
\end{equation}
for all $h\in H$ and $b\in A$. Indeed, if $\psi$ is an algebra map
then $\alpha$ is an algebra map as a composition of algebra maps.
On the other hand:
\begin{eqnarray*}
\psi (a \# h) \psi (b \# g) &=& \alpha (a) \beta (h) \alpha (b)
\beta (g)\\
\psi \bigl( (a \# h) (b \# g) \bigl ) &=& \alpha (a) \alpha
(h_{(1)} \triangleright b) \alpha\bigl(f(h_{(2)}, \,
g_{(1)})\bigl) \beta (h_{(3)}\, g_{(2)})
\end{eqnarray*}
Hence, the condition \equref{ecu111} (resp. \equref{ecu11})
follows by considering $a = b = 1_{A}$ (resp. $a = 1_A$ and $g =
1_H$) in the identity $\psi (a \# h) \, \psi (b \# g) = \psi
\bigl( (a \# h) (b \# g) \bigl )$. The converse is obvious.

Now, we prove that $\alpha : A \to A' \#^{\triangleright'}_{f'}
H'$, $\alpha (a) = u(a_{(1)}) \, \#' \, p(a_{(2)})$ is an algebra
map if and only if $p: A \to H'$ is an algebra map and (CP3)
holds. Indeed, $\alpha (ab) = \alpha (a) \alpha (b)$ is equivalent
to:
$$
u(a_{(1)} b_{(1)} ) \#' p(a_{(2)} b_{(2)} ) = u(a_{(1)}) \bigl(
p(a_{(2)}) \triangleright' u(b_{(1)}) \bigl)\, f'\bigl(p(a_{(3)}),
\, p(b_{(2)})\bigl) \, \#' \, p(a_{(4)}) \, p(b_{(3)})
$$
If we apply ${\rm Id} \ot \varepsilon_{H'}$ to this equation we
obtain (CP3), while if we apply $ \varepsilon_{A'} \ot {\rm Id}$
to the same equation we obtain that $p: A \to H'$ is an algebra
map, hence a morphism of Hopf algebras. The converse is obvious.

In a similar way we can show that the compatibility condition
\equref{ecu111} holds if and only if (CP4) and (CP6) hold. Indeed,
using the expressions of $\alpha$ and $\beta$ in terms of $(u, p)$
and respectively $(r, v)$, the equation \equref{ecu111} is
equivalent to:
\begin{eqnarray*}
&& r(h_{(1)}) \bigl(v(h_{(2)}) \triangleright' r(g_{(1)})\bigl) \,
f'\bigl(v(h_{(3)}), \, v(g_{(2)})\bigl) \, \#' \, v(h_{(4)}) \,
v(g_{(3)}) = \\
&& u\bigl(f(h_{(1)}, \, g_{(1)})\bigl) \Bigl(p\bigl(f(h_{(2)}, \,
g_{(2)})\bigl) \triangleright' r(h_{(5)}g_{(5)})\Bigl) \,
f'\Bigl(p\bigl(f(h_{(3)}, \, g_{(3)})\bigl), \,
v(h_{(6)}g_{(6)})\Bigl)  \,\, \#' \\
&& p\bigl(f(h_{(4)}, \, g_{(4)})\bigl) v(h_{(7)}g_{(7)})
\end{eqnarray*}
If  we apply ${\rm Id} \ot \varepsilon_{H'}$ to the above identity
we obtain (CP6) while if we apply $ \varepsilon_{A'} \ot {\rm Id}$
to it we get (CP4). Conversely, the compatibility condition
\equref{ecu111} follows straightforward from (CP6) and (CP4).

Finally, we prove that the commutativity condition \equref{ecu11}
is equivalent to (CP5) and (CP7). Indeed, \equref{ecu11} is
equivalent to:
\begin{eqnarray*}
&& r(h_{(1)}) \, \bigl(v(h_{(2)}) \triangleright' u(a_{(1)})\bigl)
f' \bigl(v(h_{(3)}), \, p(a_{(2)})\bigl) \#' v(h_{(4)}) \,
p(a_{(3)})
= \\
&& u(h_{(1)} \triangleright a_{(1)}) \bigl(p(h_{(2)}
\triangleright a_{(2)}) \triangleright' r(h_{(5)})\bigl) f'
\bigl(p(h_{(3)} \triangleright a_{(3)}), \,
v(h_{(6)})\bigl) \, \#' \\
&& p(h_{(4)} \triangleright a_{(4)}) \, v(h_{(7)})
\end{eqnarray*}
If  we apply ${\rm Id} \ot \varepsilon_{H'}$ to the above identity
we obtain (CP7) while if we apply $ \varepsilon_{A'} \ot {\rm Id}$
to it we get (CP5). Conversely, the commutativity condition
\equref{ecu11} follows straightforward from (CP7) and (CP5).

To conclude, we have proved that any bialgebra map $\psi: A
\#^{\triangleright}_{f} H \to A' \#^{\triangleright'}_{f'} H' $ is
uniquely determined by a quadruple $(u, p, r, v)$,  where $p: A
\to H'$ is a morphism of Hopf algebras, $u: A \to A'$, $r: H
\rightarrow A'$ and $v: H \rightarrow H'$ are unitary morphisms of
coalgebras satisfying the compatibility conditions (CP1)-(CP7)
such that $\psi: A \#^{\triangleright}_{f} H \to A'
\#^{\triangleright'}_{f'} H' $ is given by \equref{morfcros} and
the proof is finished.
\end{proof}

The compatibility conditions of \thref{toatemorf} are rather
difficult to deal with. However, there are several special cases
in which the two compatibilities simplify considerably. The first
one, which will be used in \seref{exclasros}, is the following:

\begin{corollary} \colabel{cazspecialusor}
Let $(A, H, \triangleright, f)$ and $(A, H, \triangleright', f')$
be two crossed systems of Hopf algebras such that the only Hopf
algebra map $p: A \to H$ is the trivial one. Then there exists a
bijective correspondence between the set of all morphisms of Hopf
algebras $\psi : A \#^{\triangleright}_{f} \, H \to A
\#^{\triangleright'}_{f'} \, H $ and the set of all triples $(u,
r, v)$, where $u: A \to A$, $v: H \rightarrow H$ are morphisms of
Hopf algebras, $r: H \rightarrow A$ is a unitary coalgebra map
satisfying the following compatibility conditions:
\begin{eqnarray}
r(h_{(1)}) \ot v(h_{(2)}) &{=}& r(h_{(2)}) \ot
v(h_{(1)})\eqlabel{cc1t}\\
r(h_{(1)}) \bigl(v(h_{(2)}) \triangleright' r(g_{(1)})\bigl) f'
\bigl(v(h_{(3)}), \, v(g_{(2)})\bigl) &{=}& u \bigl(f(h_{(1)}, \,
g_{(1)}) \bigl) \, r(h_{(2)} g_{(2)}) \eqlabel{cc2t}\\
r(h_{(1)}) \bigl(v(h_{(2)}) \triangleright' u(a) ) &{=}& u
(h_{(1)} \triangleright a) \, r (h_{(2)}) \eqlabel{cc3t}
\end{eqnarray}
for all $a \in A$, $g$, $h \in H$. Under the above correspondence
the morphism of Hopf algebras $\psi: A \#^{\triangleright}_{f} H
\to A \#^{\triangleright'}_{f'}  H $ corresponding to $(u, r, v)$
is given by:
\begin{equation}\eqlabel{ccmorfcrostr}
\psi(a \# h) = u(a) \, r(h_{(1)}) \, \#' \, v(h_{(2)})
\end{equation}
for all $a \in A$ and $h\in H$.

Furthermore, $\psi: A \#^{\triangleright}_{f} H \to A
\#^{\triangleright'}_{f'}  H $ given by \equref{ccmorfcrostr} is
an isomorphism if and only if $u$ and $v$ are automorphisms of
Hopf algebras.
\end{corollary}

\begin{proof}
The first part follows from \thref{toatemorf} applied for $A' =
A$, $H' = H$ and $p (a) = \varepsilon_A (a)1_H$, for all $a\in A$.
Assume now that $\psi$ is an isomorphism and let $\psi^{-1}$ be
its inverse associated to a triple $(u', r', v')$; that is
$\psi^{-1} (a \# h)
 = u'(a) \, r'(h_{(1)}) \, \#' \, v'(h_{(2)})$. Then, we have:
\begin{eqnarray}
a \# 1 &=& (\psi^{-1} \circ \psi) (a \# 1) = u'\big( u(a) \big) \#
1\eqlabel{izo1.1}\\
1 \# h &=& (\psi^{-1} \circ \psi) (1 \# h) = u' \big( r(h_{(1)})
\big) r' \big( v(h_{(2)}) \big) \# v' \big( v(h_{(3)})
\big)\eqlabel{izo1.2}
\end{eqnarray}
for all $a \in A$ and $h \in H$. Applying $Id_{A} \ot \varepsilon$
in \equref{izo1.1} (resp. $\varepsilon_{A} \otimes \Id_{H}$ in
\equref{izo1.2}), that $u' \circ u = Id_{A}$ (resp. $v' \circ v =
\Id_{H}$). In a similar manner it follows from $\psi \circ
\psi^{-1} = \Id_{A \#' H}$ that $u \circ u' = Id_{A}$ and $v \circ
v' = \Id_{H}$. Thus, $u$ and $v$ are isomorphisms.

Conversely, if $u$ and $v$ are automorphisms, then $\psi$ given by
\equref{ccmorfcrostr} is bijective with the inverse given by:
$$
\varphi : A \#' H \to A \# H, \,\,\,\, \varphi (a \#' h) =
u^{-1}(a) \big( u^{-1} \circ S \circ r \circ v^{-1} \big)
(h_{(1)}) \, \# \, v^{-1}(h_{(2)})
$$
where $u^{-1}$ (resp. $v^{-1}$) is the composition inverse of $u$
(resp. $v$).
\end{proof}

\begin{corollary} \colabel{cazspecialimport}
Let $\psi: A \#^{\triangleright}_{f} H \to A
\#^{\triangleright'}_{f'}  H $ be a Hopf algebras map between two
crossed products. Then $\psi$ stabilizes $A$ (resp. co-stabilizes
$H$) if and only if $\psi$ is a left $A$-linear (reps. right
$H$-colinear) map.
\end{corollary}

\begin{proof} A morphism $\psi : A \#^{\triangleright}_{f} \, H \to A
\#^{\triangleright'}_{f'} \, H$ stabilizes $A$ if and only if
$\psi (a \# 1_H) = a \#' 1_H$, for all $a\in A$. Taking into
account the formula for $\psi$ given by \equref{morfcros}, we
obtain that $\psi$ stabilizes $A$ if and only if $ u(a_{(1)}) \#'
p (a_{(2)}) = a \#' 1_H$, i.e. $u (a) = a$ and $p(a) =
\varepsilon_A (a) 1_{H'}$, for any $a\in A$. Thus, $\psi$ take the
form $\psi(a \# h) = a \, r(h_{(1)}) \, \#' \, v(h_{(2)})$, for
all $a \in A$ and $h\in H$ and such a map is obviously left
$A$-linear. Conversely, if $\psi$ is left $A$-linear, then $ \psi
(a \# 1) = \psi ( a \cdot (1_A \# 1_H)) = a \cdot \psi (1_A \#
1_H) = a\#' 1_H$, i.e. $\psi$ stabilizes $A$. The other statement
can be proved in a dual manner.
\end{proof}

\begin{remarks} \relabel{renouGclas}
1. Using \prref{descriereclef} and \coref{cazspecialimport} the
classification in the sense of \deref{coalspex} of all coalgebra
split extensions of $A$ by $H$ reduces to the classification of
all crossed products $A \#_{f}^{\triangleright} \, H$. The
classifying object is denoted by ${\mathcal H}^{2} (H, A)$ and was
constructed in \cite[Proposition 2.2]{agorecia} as a special case
of \cite[Theorem 3.4]{am1} as follows: let ${\mathcal C} {\mathcal
S} (A, H)$ be the set of all pairs $(\triangleright, f)$ such that
$(A, H, \triangleright, f)$ is a crossed system of Hopf algebras.
Two pairs $(\triangleright, f)$ and $(\triangleright', f') \in
{\mathcal C}{\mathcal S} (A, H)$ are called \emph{cohomologous}
and we denote this by $(\triangleright, f) \approx
(\triangleright', f')$ if there exists an unitary cocentral map
$r: H \rightarrow A$ such that
\begin{eqnarray*}
h \triangleright' a &=& r (h_{(1)}) \, (h_{(2)}
\triangleright a) \, (S_A \circ r )(h_{(3)}) \\
f'(h, \, g) &=& r(h_{(1)}) \, \bigl(h_{(2)} \triangleright
r(g_{(1)}) \bigl) \, f(h_{(3)}, \, g_{(2)})  \, (S_A \circ r)
(h_{(4)} g_{(3)})
\end{eqnarray*}
for all $a\in A$ and $h$, $g \in H$. Then \cite[Theorem 3.4]{am1}
proves that $(\triangleright, f) \approx (\triangleright', f')$ if
and only if there exists a Hopf algebra isomorphism $A
\#_{f}^{\triangleright} \, H \cong A \#_{f'}^{\triangleright'} \,
H$ that stabilizes $A$ and co-stabilizes $H$. Thus, $\approx$ is
an equivalence relation on the set ${\mathcal C} {\mathcal S} (A,
H)$. The cohomological object $ {\mathcal H}^{2} (H, A)$ is the
pointed quotient set defined by
$$
{\mathcal H}^{2} (H, A) := {\mathcal C} {\mathcal S} (A, H)/
\approx
$$

2. A more general form of classification is the following: two
coalgebra split extensions $(E, \pi)$ and $(E', \pi')$ of $A$ by
$H$ are called \emph{isomorphic} if there exists an isomorphism of
Hopf algebras $E \cong E'$. Using again \prref{descriereclef} this
classification reduces also to classifying up to a Hopf algebra
isomorphism all crossed products $A \#_{f}^{\triangleright} \, H$.
We denote by $\textsc{C}\textsc{r}\textsc{p} (H, A)$ the set of
types of Hopf algebra isomorphisms of all crossed products $A
\#_{f}^{\triangleright} \, H$ associated to all crossed systems
$(A, H, \triangleright, f)$. It is obvious that two equivalent
extensions are isomorphic and hence there exists a canonical
surjection ${\mathcal H}^{2} (H, A) \twoheadrightarrow
\textsc{C}\textsc{r}\textsc{p} (H, A)$.
\end{remarks}

\section{Classifying crossed products with the Sweedler's Hopf algebra}
\selabel{exclasros}

This section is devoted to the classification of all coalgebra
split extensions of $A$ by $H_4$, where $A$ is a Hopf algebra and
$H_4$ is the Sweedler's $4$-dimensional Hopf algebra. There are
three steps that we have to go through. First of all we have to
compute the set of all crossed systems $(A, H_4, \triangleright,
f)$ between $A$ and $H_4$. This is the computational part of our
approach. Then we describe by generators and relations all crossed
products $A\#^{\triangleright}_{f} H_4$ associated to these
crossed systems. Finally, using \thref{toatemorf}, we shall
classify the above crossed products by computing the classifying
objects ${\mathcal H}^{2} (H_4, A)$ and
$\textsc{C}\textsc{r}\textsc{p} (H_4, A)$. As a bonus of our
approach, the group of Hopf algebra automorphisms of these crossed
products is computed.

For a Hopf algebra $A$, $G(A)$ is the set of group-like elements
of $A$ and for $g$, $h \in G(A)$ we denote by $P_{g, h} (A)$ the
set of all $(g, h)$-primitive elements, that is
$$
P_{g, \, h} (A) = \{ x\in A \, | \, \Delta_A (x) = x \ot g + h \ot
x\}
$$
We denote by $P(A) = P_{1, 1} (A)$ the set of all primitive
elements of $A$ and by ${\mathcal Z} {\mathcal P} (A) := P(A) \cap
Z(A)$, where $Z(A)$ is the center of $A$.

Let $k$ be a field of characteristic $\neq 2$ and $H_{4}$ the
Sweedler's $4$-dimensional Hopf algebra having $\{1, \, g, \, x,
\, gx \}$ as a basis with the multiplication:
$$
g^{2} = 1, \quad x^{2} = 0, \quad x g = -g x
$$
and the coalgebra structure such that $g$ is a group-like element
and $x$ is $(1, g)$-primitive.

Our first classification result shows that for an arbitrary Hopf
algebra $A$ the set ${\mathcal Z} {\mathcal P} (A)$ of central
primitive elements of $A$ parameterizes all crossed systems $(A,
H_4, \triangleright, f)$.

\begin{theorem} \thlabel{ex1thgen}
Let $k$ be a field of characteristic $\neq 2$ and $A$ a Hopf
algebra. Then there exists a bijection between the set of all
crossed systems $(A, H_4, \triangleright, f)$ and the set
${\mathcal Z} {\mathcal P} (A)$ of all central primitive elements
of $A$ such that the crossed system $(A, H_4, \triangleright, f)$
corresponding to $a \in {\mathcal Z} {\mathcal P} (A)$ is given as
follows: the action $ \triangleright : H_4 \ot A \to A$ is the
trivial action $h \triangleright b = \varepsilon(h) \, b$, for any
$h \in H_4$, $b \in A$ and the cocycle $f = f_a : H_4 \ot H_4 \to
A$ is given by the following formula:
\begin{equation} \eqlabel{cocicprim}
\begin{tabular}{c|cccc}
  $f$ & 1 & $g$ & $x$ & $gx$ \\
  \hline
  1     & 1 & 1   & 0   &  0 \\
  $g$   & 1 & 1   & 0   &  0 \\
  $x$   & 0 & 0   & a   & -a \\
  $gx$  & 0 & 0   & a   & -a \\
\end{tabular}
\end{equation}
In particular, if ${\mathcal Z} {\mathcal P} (A) = \{0\}$, then
there are no nontrivial crossed systems of Hopf algebras $(A, H_4,
\triangleright, f)$ and thus the only crossed product $A
\#^{\triangleright}_{f} H_4$ is the usual tensor product $A \ot
H_4$ of Hopf algebras, that is ${\mathcal H}^{2} (H, A) \cong
\textsc{C}\textsc{r}\textsc{p} (H, A) = \{ A \ot H_4 \}$.
\end{theorem}

\begin{proof} We shall compute all crossed systems
$(A, H_4, \triangleright, f)$: i.e. we have to describe all
coalgebra maps $ \triangleright : H_4 \ot A \to A$, $f : H_4 \ot
H_4 \to A$ satisfying the compatibility conditions \equref{1a} -
\equref{7}. Let $(\triangleright, f)$ be such a pair. We shall
prove first that $ \triangleright : H_4 \ot A \to A$ is
necessarily the trivial action. Indeed, let $a\in A$. If we apply
the compatibility condition \equref{6} for $g := x$ we obtain,
taking into account \equref{1a},
$$
x \otimes a + g \otimes ( x\triangleright a) = 1 \otimes (x
\triangleright a) + x \otimes (g \triangleright a)
$$
If we apply $x^* \ot {\rm Id}_A$ and $g^* \ot {\rm Id}_A $ to this
equation (where $x^*$ and $g^* \in H_4^*$ are the elements of the
dual basis of $\{1, \, g, \, x, \, gx \}$) we obtain that $g
\triangleright a = a$ and $x \triangleright a = 0$. Now, if we
apply \equref{6} for $g := gx$ we obtain, using $g \triangleright
a = a$, that $1 \otimes (gx) \triangleright a = g \otimes (gx)
\triangleright a $ and hence $(gx) \triangleright a = 0$. Thus we
have proved that $\triangleright: H_4 \ot A \to A $ acts trivially
and hence the compatibility conditions \equref{1a} - \equref{2}
and \equref{6} are trivially fulfilled.

It remains to describe all the cocycles $f: H_4 \ot H_4 \to A$.
Since the action $\triangleright$ is trivial, the compatibility
condition \equref{4} takes the form, $b \, f (h, h') = f(h, h') \,
b$, for all $b\in A$, $h$, $h' \in H_4$, i.e. \equref{4} is
equivalent to the fact that ${\rm Im} (f) \subseteq Z(A)$.
Furthermore, the normalizing condition \equref{3} is equivalent to
$$
f(1,1) = f (1, g) = f (g, 1) = 1, \quad f (1, x) = f (1, gx) = f
(x, 1) = f( gx, 1) = 0
$$
From now on we assume that $f$ is such a normalized map. The next
step proves that the compatibility condition \equref{7} holds for
$f$ if and only if
\begin{equation} \eqlabel{cocih4tri}
f(g, g) = 1, \qquad f(x, g) = f(g, x) = f(gx, g) = f(g, gx) = 0
\end{equation}
Indeed, first of all we observe that \equref{7} is trivially
fulfilled for $g = 1$ or $h = 1$. Now, the compatibility condition
\equref{7} holds for $(x, g)$, $(gx, g)$, $(g, x)$ and
respectively $(g, gx)$ if and only if:
\begin{eqnarray*}
xg \otimes 1 + 1 \otimes f(x, g) &=& g \otimes f(x, g) + xg
\otimes f(g, g) \\
-x \otimes f(g, g) + g \otimes f(gx, g) &=& 1 \otimes f(gx, g) -
x \otimes 1 \\
gx \otimes 1 + 1 \otimes f(g, x) &=& g \otimes f(g,x) + gx \otimes
f(g,g)\\
x \otimes f(g, g) + g \otimes f(g, gx) &=& 1 \otimes f(g, gx) + x
\otimes 1
\end{eqnarray*}
These four equations are equivalent to the fact that
\equref{cocih4tri} holds. Now, with the values of $f$ given by
\equref{cocih4tri}, it is just a straightforward computation to
prove that \equref{7} holds for any other pair of $\{g, x, gx \}
\times \{g, x, gx \}$. More precisely, there are other five
possibilities that need to be checked, namely $(g, g)$, $(x, x)$,
$(x, gx)$, $(gx, x)$ and $(gx, gx)$. In all these cases we obtain
compatibilities that are trivially fulfilled: in $(x, x)$ we
obtain $1 \ot f (x, x) = 1 \ot f (x, x)$, in $(x, gx)$ we obtain $
g \ot f (x, gx) = g \ot f (x, gx)$, in $(gx, x)$ we obtain $ g \ot
f (gx, x) =  g \ot f (gx, x)$ and finally in $(gx, gx)$ we obtain
$ 1 \ot f (gx, gx) =  1 \ot f (gx, gx)$.

To summarize, we have proved so far that a pair of maps
$(\triangleright, f)$ satisfies all axioms of a crossed system
\equref{1a} - \equref{7}, except for the cocycle axiom \equref{5},
if and only if $\triangleright: H_4 \ot A \to A$ is the trivial
action and $f: H_4 \ot H_4 \to A$ is given by
$$
\begin{tabular}{c|cccc}
  $f$ & 1 & $g$ & $x$ & $gx$ \\
  \hline
  1     & 1 & 1   & 0   & 0 \\
  $g$   & 1 & 1   & 0   & 0 \\
  $x$   & 0 & 0   & ?   & ? \\
  $gx$  & 0 & 0   & ?   & ? \\
\end{tabular}
$$

In the final step of our investigation the above values of $f$
marked with $?$ will be determined such that the cocycle condition
\equref{5} holds. First of all, using the values of $f$ that we
have already obtained and the fact that $f$ is a morphism of
coalgebras one can prove that $f(x, x)$, $f(x, gx)$, $f(gx, x)$
and $f(gx, gx)$ are primitive elements of $A$. Indeed, if we write
down the condition $\Delta_A (f (x, x)) = (f \ot f ) \circ
\Delta_{H_4 \ot H_4} (x \ot x)$ we obtain that:
\begin{eqnarray*}
\Delta_A (f (x, x)) &=& f(x, x) \ot 1 + f(x, g) \ot f(1, x) + f(g,
x) \ot f(x, 1) + f (g, g) \ot f(x, x) \\
&=& f(x, x) \ot 1 +  1 \ot f(x, x)
\end{eqnarray*}
that is $f(x, x)$ is a primitive element of $A$. In the same way
we can prove that $f(x, gx)$, $f(gx, x)$ and $f(gx, gx)$ are also
primitive elements.

The only compatibility which remains to be fulfilled by $f$ is the
cocycle condition \equref{5} which, considering that the action
$\triangleright$ is trivial, takes the simplified form:
\begin{equation}\eqlabel{5simpli}
f(h_{(1)}, \, l_{(1)}) \, f\bigl(y, \, h_{(2)} l_{(2)} \bigl) =
f(y_{(1)}, \, h_{(1)}) \, f(y_{(2)}h_{(2)}, \, l)
\end{equation}
for all $h$, $l$, $y\in H_4$. First we observe that the condition
\equref{5simpli} holds if one of the elements $h$, $l$ or $y$ is
$1$. Thus, \equref{5simpli} holds if and only if it holds in all
triples $(h, l, y) \in \{g, x, gx \}^3$. Thus, there are $27$
equations that have to be fulfilled. However, this is a routinely
check so we indicate only the main steps of the proof. First we
observe that \equref{5simpli} holds for the triple $(h, l, y)$
equal to $(g, x, x)$, $(g, x, gx)$ and respectively $(x, g, x)$ if
and only if
\begin{equation} \eqlabel{teststres}
f (gx, x) = - f( x, gx), \quad f (gx, gx) = - f (x,x), \quad f(x,
gx) = - f (x, x)
\end{equation}
Thus, if we denote $f(x, x) = a \in Z(A) \cap P(A)$, we obtain
that $f( gx, x) = a$ and $f (x, gx) = f(gx, gx) = - a$ i.e.
\equref{cocicprim} holds. Now, by a long but straightforward
computation one can see that for the other $24$ possibilities of
choosing the triple $(h, l, y)$, the cocycle condition
\equref{5simpli} is either trivially fulfilled or equivalent to
one of the compatibilities \equref{teststres}. For instance, if we
consider the triple $(h, l, y) := (x, g, gx)$, \equref{5simpli}
holds if and only if  $f (gx, gx) = - f( gx, x)$, which can be
obtained from \equref{teststres}. On the other hand, using the
values of $f$ that we have already determined, it is easy to see
that for the triple $(h, l, y) := (x, x, gx)$, \equref{5simpli}
holds automatically. The last statement follows from the first
part: if $a = 0$, then $f = f_0$ is just the trivial cocycle $f
(g, h) = \varepsilon (g) \varepsilon (h) 1_A$, for all $g$, $h \in
H_4$ and hence the proof is finished.
\end{proof}

\begin{examples} \exlabel{singurulcrs}
$(1)$ Let $A$ be a finite dimensional Hopf algebra over a field of
characteristic zero. Then $P(A) = \{0\}$ \cite[Exercise
4.2.16]{dnr} and hence there are no nontrivial crossed systems of
Hopf algebras $(A, H_4, \triangleright, f)$. The conclusion fails
if ${\rm char} (k) = p \geq 3 $: in \coref{nrtip1} we shall prove
that there exists an infinite number of types of isomorphisms of
crossed products $A \# H_4$.

$(2)$ Let $A = U(L)$ be the enveloping algebra of a Lie algebra
$L$ over a field $k$ of characteristic zero. Then, $P (U (L)) = L$
\cite[Proposition 5.5.3]{M}. Thus, the set of all primitive
central elements of $U(L)$ if given by
$$
{\mathcal Z} {\mathcal P} (U(L)) = \{ l \in L \, | \,\, [l, y] =
0, \, \forall \, y \in L \} = Z (L)
$$
where $Z (L)$ is the center of the Lie algebra $L$. Typical
examples of Lie algebras with non-trivial center are the nilpotent
Lie algebras. They will provide nontrivial examples of crossed
products $U(L) \#H_4$.
\end{examples}

For any $a \in {\mathcal Z} {\mathcal P} (A)$ we shall denote by
$A_{(a)} := A \# H_4$, the crossed product associated to the
crossed system constructed in \thref{ex1thgen}. In what follows we
will describe explicitly by generators and relations the Hopf
algebras $A_{(a)}$. As the Hopf algebra map $A \to A \# H_4$, $z
\mapsto z \# 1$, is injective we shall identify $ z = z \# 1$, for
any $z \in A$. Let $\{e_i \}_{i\in I}$ be a $k$-basis of $A$ and
we denote
$$
g = 1 \# g, \quad x:= 1 \# x, \quad w := 1 \# gx \, \in A_{(a \, |
\, g, \, x)}
$$
We can easily show, by using formula \equref{cocicprim}, that $w =
g x $ in $A_{(a)}$. Indeed,
$$
g x  = (1 \# g) (1 \# x) = f (g, x_{(1)}) \# g x_{(2)} = f(g, x)
\# g + f(g, g) \# gx = 1 \# gx = w
$$
Similarly, we can prove that the following relations hold in the
Hopf algebra $A_{(a)}$:
$$
g^2 = 1, \quad x^2 = a, \quad x g = - g x, \quad g \, e_i = e_i \,
g, \quad x \, e_i = e_i \, x
$$
for any $i\in I$. Indeed, for example (below we denote $x' = x$)
we have:
\begin{eqnarray*}
x^2 &=& (1\# x) (1\# x) = f(x_{(1)}, x_{(1)}') \# x_{(2)} x_{(2)}'
\\
&=& f(x, x) \#1 + f(x, g) \# x + f(g, x) \# x + f(g, g) \# x^2\\
&=& a \#1 = a
\end{eqnarray*}
and
$$
x e_i = (1 \# x) (e_i \# 1) = e_i \# x = (e_i \# 1) (1 \# x) = e_i
x
$$
Now, as a vector space, $A_{(a)} := A \ot H_4$, hence, the set
$\{e_i, \, e_i g, \, e_i x, \, e_i gx \, | \, i\in I \}$ is a
$k$-basis of $A_{(a)}$, where we identify $e_i = e_i \# 1$, $e_i
\# g = (e_i \# 1) (1 \# g) = e_i g$ and so on.

Using \thref{ex1thgen} and the above computations we obtain:

\begin{corollary} \colabel{noucorH}
Let $k$ be a field of characteristic $\neq 2$, $A$ a Hopf algebra
and consider $\{e_i \}_{i\in I}$ a $k$-basis of $A$. Then a Hopf
algebra $E$ is isomorphic to a crossed product of Hopf algebras $A
\# H_4$ if and only if $E \cong A_{(a)}$, for some $a \in
{\mathcal Z} {\mathcal P} (A)$, where $A_{(a)}$ is the Hopf
algebra having $\{e_i, \, e_i g, \, e_i x, \, e_i gx \, | \, i\in
I \}$ as a $k$-basis and the multiplication is subject to the
following relations for any $i\in I$:
\begin{equation} \eqlabel{relatiidinou}
g^2 = 1, \quad x^2 = a, \quad x g = - g x, \quad g \, e_i = e_i \,
g, \quad x \, e_i = e_i \, x
\end{equation}
while the multiplications in $A_{(a)}$ of two elements $e_i$ and
$e_j$ is the same to the one in $A$. The coalgebra structure and
the antipode of $A_{(a)}$ are given by:
$$
\Delta (e_i) := \Delta_A (e_i), \quad \Delta (g) := g\ot g, \quad
\Delta (x) := x \ot 1 + g \ot x
$$
$$
\varepsilon (e_i) := \varepsilon_A (e_i), \quad \varepsilon (g) :=
1, \quad \varepsilon (x) := 0
$$
$$
S (e_i) := S_A (e_i), \quad  S (g) := g, \quad S (x) := - g x
$$
for all $i\in I$.
\end{corollary}

\begin{example} \exlabel{lie}
Let $L := h(n, k)$ be the $2n+1$ dimensional Heisenberg Lie
algebra over a field $k$ of characteristic zero. That is, $h(n,
k)$ has a basis $\{x_1, \cdots, x_n, y_1, \cdots, y_n, z \}$ and
the only non-zero Lie brackets are $[x_i, y_i] = z$, for all $i =
1, \cdots, n$. Since ${\rm char} (k) = 0 $, it is well known that
$Z ( h(n, k) ) = \{ q z \, | \, q \in k \}$. Let $A := U ( h(n,
k))$. Then ${\mathcal Z} {\mathcal P} (A) = k \, z $ and thus any
primitive central element $ l \in U ( h(n, k))$ is of the form $l
= q z$, for some scalar $q \in k$. If $q = 0$, then $U( h(n, k)
)_{(0)} = U( h(n, k) ) \ot H_4$, the tensor product of the two
Hopf algebras. Assume now that $q \neq 0$. Using \coref{noucorH}
we can easily prove that $U( h(n, k) )_{( q z)} = k_q \lan x_i, \,
y_i, \, g, \, x \, | \, i = 1, \cdots, n \ran$, where by $k_q \lan
x_i, \, y_i, \, g, \, x \, | \, i = 1, \cdots, n \ran$ we denote
the quantum group generated as an algebra by $\{ x_i, \, y_i, \,
g, \, x \, | \, i = 1, \cdots, n \}$ subject to the following
relations for any $i = 1, \cdots, n$:
$$
g^2 = 1, \quad xg = -gx, \quad gx_i = x_i g, \quad g y_i = y_i g
$$
$$
x x_i = x_i x, \quad xy_i = y_i x, \quad x_i y_i - y_i x_i =
q^{-1} \, x^2
$$
with the coalgebra structure given such that $g$ is a group-like
element, $x$ is a $(1, g)$-primitive element and $x_i$, $y_i$ are
primitive elements for any $i = 1, \cdots, n$. We observe that
there exists a Hopf algebra isomorphism $k_q \lan x_i, \, y_i, \,
g, \, x \, | \, i = 1, \cdots, n \ran \cong k_1 \lan x_i, \, y_i,
\, g, \, x \, | \, i = 1, \cdots, n \ran$, for any $q \in k^*$ and
hence $\textsc{C}\textsc{r}\textsc{p} ( \, U ( h(n, k)), \, H_4)$
has two elements.
\end{example}

Now, we shall classify $A_{(a)}$ for those Hopf algebras $A$ such
that the only Hopf algebra map $f : A \to H_4$ is the trivial one,
namely $f (z) = \varepsilon (z) 1$, for all $z \in A$. The typical
example is again $A = U (L)$, for a Lie algebra $L$. We recall
from \cite[Lemma 4.6]{ABM1} that $\Aut_{\rm Hopf}(H_4) \cong k^*$:
explicitly, any automorphism $v : H_4 \to H_4$ is of the form
\begin{equation} \eqlabel{impleh4}
v(g) = g, \quad v (x) = \beta \, x, \quad v(gx) = \beta \, gx
\end{equation}
for some non-zero scalar $\beta \in k^*$. It what follows, the
automorphism $v$ of $H_4$ implemented by $\beta \in k^*$ as in
\equref{impleh4} will be denoted by $v = v_{\beta}$.

\begin{theorem} \thlabel{clssifica1ho}
Let $k$ be a field of characteristic $\neq 2$ and $A$ a Hopf
algebra such that the only Hopf algebra map $A \to H_4$ is the
trivial one. Let $a$, $b\in {\mathcal Z} {\mathcal P} (A)$ be two
central primitive elements of $A$. Then there exists a bijection
between the set of all Hopf algebra isomorphisms $\psi : A_{(a)}
\to A_{(b)} $ and the set of all pairs $(u, \, \beta) \in
\Aut_{\rm Hopf}(A) \times k^*$ such that $u(a) = \beta^{2} \, b$.
\end{theorem}

\begin{proof} We show that the bijection
from the statement is given such that the isomorphism $\psi =
\psi_{u, \beta} : A_{(a)} \to A_{(b)}$ corresponding to a pair
$(u, \beta) \in \Aut_{\rm Hopf}(A) \times k^*$ satisfying $u(a) =
\beta^{2} b$, is given by
\begin{equation} \eqlabel{293}
\psi_{u, \beta} (z \# h) = u(z) \# v_{\beta}(h)
\end{equation}
for all $z \in A$, $h \in H_{4}$, where $v_{\beta} : H_4 \to H_4$
is the automorphism given by \equref{impleh4}.

Indeed, $A_{(a)} = A \#_a \, H_4$, where by $A \#_a \, H_4$ we
denoted the crossed product associated to the central primitive
element $a$ from \thref{ex1thgen}. It follows from
\coref{cazspecialusor} that the set of all Hopf algebra maps $\psi
: A_{(a)} \to A_{(b)} $ is in bijection with the set of all
triples $(u, r, v)$, where $u: A \to A$, $v : H_{4} \to H_{4}$ are
Hopf algebra maps and $r: H_{4} \to A$ is a unitary coalgebra map
such that the compatibility conditions \equref{cc1t}-\equref{cc3t}
are fulfilled. Under this bijection the morphism $\psi = \psi_{(u,
r, v)}$ corresponding to $(u, r, v)$ is given by
\equref{ccmorfcrostr}. We will prove now that, under this
bijection, the isomorphisms $\psi = \psi_{(u, r, v)}$ correspond
precisely to the triples $(u, r, v)$ such that $u: A \to A$, $v:
H_4 \to H_4$ are isomorphisms of Hopf algebras, $r: H_4 \to A$ is
the trivial coalgebra map and $u(a) = \beta^{2} b$, where $\beta
\in k^*$ is the scalar that implements the isomorphism $v =
v_{\beta}$ of $H_4$.

Indeed, it follows from \coref{cazspecialusor} that $\psi =
\psi_{(u, r, v)}$ is an isomorphism if and only if $u: A \to A$
and $v: H_4 \to H_4$ are isomorphisms of Hopf algebras. Let $(u,
v) \in \Aut_{\rm Hopf}(A) \times \Aut_{\rm Hopf}(H_4)$ and $\beta
\in k^*$ such that $v = v_{\beta}$. We will show now that the
compatibility conditions \equref{cc1t}-\equref{cc3t} are fulfilled
for the triple $(u, r, v)$ if and only if $r: H_4 \to A$ is the
trivial coalgebra map, i.e. $r(h) = \varepsilon (h) 1_A$, for all
$h\in H_4$ and $u(a) = \beta^{2} b$ and this will finish the
proof. First we remark that if $r : H_4 \to A$ is the trivial
coalgebra map then the compatibility condition \equref{cc1t}
holds. Conversely, since $v$ is an isomorphism the compatibility
condition \equref{cc1t} is equivalent to
\begin{equation} \eqlabel{cc1ta}
r(h_{(1)}) \ot h_{(2)} = r(h_{(2)}) \ot h_{(1)}
\end{equation}
for any $h\in H_4$, i.e. to the fact that $r: H_4 \to A$ is a
cocentral map. Applying \equref{cc1ta} for $h = x$ and then for $h
= gx$ we obtain
\begin{eqnarray*}
r(x) \ot 1 + r(g) \ot x &=& 1 \ot x + r(x) \ot g\\
r(gx) \ot g + 1 \ot gx &=& r(g) \ot gx + r(gx) \ot 1
\end{eqnarray*}
It follows from here that $r(x) = 0$, $r(g) = 1$ and $r(gx) = 0$.
Thus $r$ is the trivial map. Thus \equref{cc1t} holds for a triple
$(r, u, v)$ such that $\psi = \psi_{(u, r, v)}$ is an isomorphism
if and only if $r$ is a trivial map. In this case, the
compatibility condition \equref{cc3t} also holds since $r: H_4 \to
A$ is the trivial map and the actions $\triangleright$ and
$\triangleright'$ of the crossed systems $ (A, H_4,
\triangleright, f = f_a)$ and $ (A, H_4, \triangleright', f = f_b
)$ are also the trivial actions according to \thref{ex1thgen}.

Now, we look at the remaining compatibility condition
\equref{cc2t} that has to be satisfied by a triple $(u, r, v)$ for
which $\psi_{(u, r, v)}$ is an isomorphism. Since $r$ is the
trivial coalgebra map and $\triangleright$ and $\triangleright'$
are the trivial actions, the compatibility condition \equref{cc2t}
is equivalent to
$$
f_{b} \big( v_{\beta} (h), v_{\beta} (k) \big) = u \big( f_{a}(h,
k) \big)
$$
for all $h$, $k \in H_{4}$. It is easy to see that this
compatibility condition holds if and only if $\beta^{2} b = u(a)$.
For example, $f_{b} \big( v_{\beta}(x), v_{\beta}(x) \big) =
\beta^{2} f_{b}(x, x) = \beta^{2}b$ and $u \big( f_{a}(x,x) \big)
= u(a)$, so \equref{cc2t} holds for $(h, k) = (x, x)$ if and only
if $\beta^{2} b = u(a)$. The rest is straightforward and the proof
is now complete.
\end{proof}

\begin{remark} \relabel{calculefec}
Let $A$ be a Hopf algebra such that the only Hopf algebra map $A
\to H_4$ is the trivial one. For any $a\in {\mathcal Z} {\mathcal
P} (A)$ we define
$$
G(a) := \{ \, (u, \, \beta) \in \Aut_{\rm Hopf}(A) \times k^* \,
\, | \, \,  u(a) = \beta^{2} \, a \}
$$
Then $G(a)$ is a subgroup of $\Aut_{\rm Hopf}(A) \times k^*$ and
it follows from \thref{clssifica1ho} that there exists an
isomorphism of groups
$$
G(a) \,\, \cong \,\, \Aut_{\rm Hopf}(A_{(a)}), \qquad (u, \beta)
\mapsto \psi_{u, \beta}
$$
where $\psi_{u, \beta}$ is given by \equref{293}.
\end{remark}

Now, summarizing our results, we obtain the following
classification result:

\begin{theorem} \thlabel{classinteza1}
Let $k$ be a field of characteristic $\neq 2$ and $A$ a Hopf
algebra such that the only Hopf algebra morphism $A \to H_4$ is
the trivial one. Then:

$(1)$ There exists a bijection between the set of types of Hopf
algebra isomorphisms of all crossed products $A \#H_4$ and the
quotient pointed set ${\mathcal Z} {\mathcal P} (A)/\equiv $,
where $ \equiv$ is the equivalence relation on ${\mathcal Z}
{\mathcal P} (A)$ defined by: $a \equiv b$ if and ond only if
there exists a pair $(u, \, \beta) \in \Aut_{\rm Hopf}(A) \times
k^*$ such that $u(a) = \beta^{2} \, b$. Hence,
$\textsc{C}\textsc{r}\textsc{p} (H_4, A) \cong {\mathcal Z}
{\mathcal P} (A)/\equiv $.

$(2)$ There exists a bijection ${\mathcal H}^{2} (H_4, A) \cong
{\mathcal Z} {\mathcal P} (A)$.
\end{theorem}

\begin{proof}
It follows from \thref{clssifica1ho} and \thref{ex1thgen} once we
observe that the isomorphism $\psi_{u, \beta}$ given by
\equref{293} stabilizes $A$ (resp. co-stabilizes $H_4$) if and
only if $u = {\rm Id}_A$ (resp. $\beta = 1$), i.e. two coalgebra
split extensions of $A$ by $H_4$, $(A_{(a)}, \, \pi_{H_4})$ and
$(A_{(b)}, \, \pi_{H_4})$ are equivalent if and only if $b = a$.
\end{proof}

\section{Examples} \selabel{exdetaliate}
In this section we provide some explicit classification for all
crossed products of the form $A \, \#^{\triangleright}_{f} H_{4}$,
where $A$ is the polynomial Hopf algebra $k[Y]$ ($Y$ is a
primitive element) or two of its quotients when ${\rm char} (k) =
p \geq 3$: the $p$-dimensional Hopf algebras $k \lan y ~|~ y^{p} =
0 \ran$ and $k \lan y ~|~ y^{p} = y \ran$. To start with, we
collect some technical results. First, if $A$ is a Hopf algebra
generated by primitive elements then the only Hopf algebra
morphism $p: A \to H_{4}$ is the trivial one, i.e. $p(a) =
\varepsilon(a) 1$, for all $a \in A$, since $P(H_{4}) = \{0\}$. We
also need to know the set of primitive elements and the
automorphism groups of these Hopf algebras. For $k[Y]$ the
description of the primitive elements follows from a more general
result for universal enveloping algebras \cite[Proposition
5.5.3]{M}: $P \big (k[Y] \big) = \{q \, Y ~|~ q \in k \}$ if ${\rm
char} (k) = 0$ and $P \big( k[Y] \big)$ is the span of all
$Y^{p^{i}}$, $i \geq 0$, if ${\rm char} (k) = p \geq 2$. For the
description of the primitive elements of the other Hopf algebras
and of the groups of automorphisms we have the following:

\begin{lemma}\lelabel{lemacumorfisme}
Let $k$ be a field of characteristic $\neq 2$. Then:

$(1)$ If ${\rm char} (k) = p \geq 3$ then $P\bigl(k \langle y ~|~
y^{p} = 0 \rangle \bigl) = \{q \, y ~|~ q \in k\}$ and $P\bigl(k
\langle y ~|~ y^{p} = y \rangle \bigl) = \{q \, y ~|~ q \in k\}$;

$(2)$ $u: k[Y] \to k[Y]$ is a Hopf algebra automorphism if and
only if there exists an $\alpha \in k^{*}$ such that $u(Y^{i}) =
\alpha^{i} \, Y^{i}$ for all $i \geq 0$. In particular, $\Aut_{\rm
Hopf}(k[Y]) \simeq k^{*}$;

$(3)$ If ${\rm char} (k) = p \geq 3$ then $u : k \lan y ~|~ y^{p}
= 0 \ran \to k \lan y ~|~ y^{p} = 0 \ran$ is a Hopf algebra
automorphism if and only if there exists an $\alpha \in k^*$ such
that $u(y^{j}) = \alpha^{j} y^{j}$, for all $j \geq 0$. In
particular, $\Aut_{\rm Hopf} \big( k \lan y ~|~ y^{p} = 0 \ran
\big) \simeq k^*$;

$(4)$ If ${\rm char} (k) = p \geq 3$ then $u : k \lan y ~|~ y^{p}
= y \ran \to k \lan y ~|~ y^{p} = y \ran$ is a Hopf algebra
automorphism if and only if there exists an $\alpha \in
\mathbb{F}_{p}^*$ such that $u(y^{j}) = \alpha^{j} y^{j}$, for all
$j \geq 0$. In particular, $\Aut_{\rm Hopf} \big( k \lan y ~|~
y^{p} = y \ran \big) \simeq \mathbb{F}_{p}^*$, where
$\mathbb{F}_{p}$ is the field with $p$ elements.
\end{lemma}

\begin{proof} The proof is a straightforward
computation. We only remark for $(4)$ that if $u : k \lan y ~|~
y^{p} = y \ran \to k \lan y ~|~ y^{p} = y \ran$ is an automorphism
then $u(y) = \alpha y$, for some $\alpha \in k^*$ that must
satisfy $\alpha^{p} = \alpha$, since
$$
\alpha y = u(y) = u(y^{p}) = \alpha^{p} y^{p} = \alpha^{p} y
$$
As $\alpha \neq 0$, $\alpha^{p-1} = 1$. Taking into account that
$\{ \alpha \in k ~|~ \alpha^{p-1} = 1 \} = \mathbb{F}_{p}^*$, the
conclusion follows.
\end{proof}

In order to classify all crossed products $k[Y] \# H_4$ we
distinguish two cases depending on the characteristic of the base
field $k$.

\begin{proposition}\prlabel{exemplul1}
Let $k$ be a field of characteristic zero and $A := k[Y]$ the
polynomial Hopf algebra. Then:

$(1)$ Up to an isomorphism of Hopf algebras there exist exactly
two crossed products of Hopf algebras $A \# H_4$: $A \ot H_4$ and
$A_{(\infty)}$, where by $A_{(\infty)}$ we denote the infinite
dimensional Hopf algebra generated by $g$ and $x$ subject to the
relations:
\begin{equation} \eqlabel{relatiiex1}
\quad g^2 = 1, \quad x g = - g x
\end{equation}
and with the coalgebra structure given such that $g$ is a
group-like element and $x$ is $(1, g)$-primitive. Furthermore, we
have the following isomorphisms of groups:
$$
\Aut_{\rm Hopf} \big( A \otimes H_{4} \big) \simeq k^* \times k^*
\qquad {\rm and} \qquad \Aut_{\rm Hopf} \big( A_{(\infty)} \big)
\simeq k^*
$$

$(2)$ There exists a bijection ${\mathcal H}^{2} (H_4, A) \cong
k$.
\end{proposition}

\begin{proof}
According to \coref{noucorH}, a Hopf algebra is isomorphic to a
crossed product of Hopf algebras $A \# H_4$ if and only if $E
\simeq A_{(qY)}$ for some $q \in k$, where $A_{(qY)}$ is the
infinite dimensional Hopf algebra generated by $Y$, $g$ and $x$,
subject to the relations
\begin{equation}
g^{2} = 1, \quad  x^{2} = qY, \quad  xg = - gx, \quad gY = Yg,
\quad xY = Yx
\end{equation}
If $q \neq 0$ then, among the previous relations, the ones in
\equref{relatiiex1} are independent. The final statement of $(1)$
follows \reref{calculefec} and \leref{lemacumorfisme}. The
classification part is a consequence of \thref{classinteza1}.
\end{proof}

The case $\mbox{char}(k) = p \geq 3$ is more interesting. We
denote by $k^{(\NN)}$ the set of sequences with finitely many
non-zero terms from $k$. For $(\alpha_{i})_{i \geq 0} \in
k^{(\NN)}$ for which there exists $i \geq 0$ such that $\alpha_{i}
\neq 0$ we define:
$$
G \big( (\alpha_{i})_{i \geq 0} \big) := \{ (\alpha, \beta) \in
k^{*} \times k^{*} \, | \, \alpha^{p^{i}} = \beta^{2}, \, {\rm
for} \, {\rm all} \, i \, {\rm such} \, {\rm that} \,\, \alpha_{i}
\neq 0 \}
$$
It is easily seen that $G \big( (\alpha_{i})_{i \geq 0} \big)$ is
a subgroup of $k^* \times k^*$.

\begin{proposition}\prlabel{exemplul2}
Let $k$ be a field of characteristic $p \geq 3$ and $A:= k[Y]$ the
polynomial Hopf algebra. Then:

$(1)$ A Hopf algebra $E$ is isomorphic to a crossed product $A \#
H_{4}$ if and only if $E \simeq A_{( (\alpha)_{i \geq 0})}$, for
some $(\alpha_{i})_{i\geq 0} \in k^{(\mathbb{N})}$, where $A_{(
(\alpha)_{i \geq 0})}$ is the infinite dimensional quantum group
generated by $Y$, $g$ and $x$, subject to the following relations
\begin{equation} \eqlabel{relatiiex2}
g^{2} = 1, \quad  x^{2} = \sum_{i \geq 0} \alpha_{i}Y^{p^{i}},
\quad xg = - gx, \quad gY = Yg, \quad xY = Yx
\end{equation}
and with the coalgebra structure given such that $Y$ is a
primitive element, $g$ is a group-like element and $x$ is $(1,
g)$-primitive. Furthermore, we have the following isomorphisms of
groups:
$$
\Aut_{\rm Hopf} \big( A \ot H_{4} \big) \simeq k^* \times k^*
\qquad {\rm and} \qquad \Aut_{\rm Hopf} \big( A_{((\alpha_{i})_{i
\geq 0})} \big) \simeq G \big( (\alpha_{i})_{i \geq 0} \big)
$$
for all $(\alpha_{i})_{i \geq 0} \in k^{(\NN)} \setminus \{0\}$ .

$(2)$ There exists a bijection $\textsc{Crp} (H_4, \, k[Y]) \cong
k^{(\mathbb{N})}/_{\sim}$, where $\sim$ is the equivalence
relation on $k^{(\mathbb{N})}$ defined by: $(\alpha_{i})_{i \geq
0} \sim (\beta_{i})_{i \geq 0}$ if and only if there exists
$(\alpha, \beta) \in k^* \times k^*$ such that $\alpha^{p^{i}}
\alpha_{i} = \beta^{2} \beta_{i}$, for all $i \geq 0$.

$(3)$ There exists a bijection ${\mathcal H}^{2} (H_4, k[Y] )
\cong k^{(\mathbb{N})}$.
\end{proposition}

\begin{proof}
The set of all primitive elements of $A = k[Y]$ is $P(A) = \{
\sum_{i \geq 0}\alpha_{i}Y^{p^{i}} ~|~ (\alpha_{i})_{i \geq 0} \in
k^{(\NN)} \}$. Denoting by $A_{( (\alpha)_{i \geq 0})}$ the Hopf
algebra associated to the primitive element $\sum_{i \geq 0}
\alpha_{i}Y^{i}$, where $(\alpha_{i})_{i \geq 0} \in
k^{(\mathbb{N})}$, we obtain the first part of $(1)$. The
description of the automorphism groups of $A_{( (\alpha_{i})_{i
\geq 0})}$ follows from \reref{calculefec} and
\leref{lemacumorfisme}; indeed, we have
$$
\Aut_{\rm Hopf} \big( A \otimes H_{4} \big) \simeq \Aut_{\rm
Hopf}(A) \times k^* \simeq k^* \times k^*
$$
If $(\alpha_{i})_{i \geq 0}$ is not the null sequence then, by
\reref{calculefec}, $A_{( (\alpha_{i})_{i \geq 0})} \simeq G \Big(
\sum_{i \geq 0} \alpha_{i}Y^{p^{i}} \Big)$, where
$$
G \Big( {\sum_{i \geq 0}\alpha_{i}Y^{p^{i}}} \Big) = \Big\{ (u,
\beta) \in \Aut_{\rm Hopf}(A) \times k^* \,\, | \,\, u \Big(
\sum_{i \geq 0} \alpha_{i} Y^{p^{i}} \Big) = \beta^{2} \sum_{i
\geq 0} \alpha_{i} Y^{p^{i}} \Big\}
$$
Taking into account the description of $\Aut_{\rm Hopf}(A)$ given
in \leref{lemacumorfisme}, we obtain
$$
G \Big({\sum_{i \geq 0}\alpha_{i}Y^{p^{i}}} \Big) \simeq \Big\{
(\alpha, \beta) \in k^* \times k^* \,\, | \,\, \alpha^{p^{i}}
\alpha_{i} = \beta^{2} \alpha_{i} \,\, {\rm for} \,\, {\rm all}
\,\, i \Big\} = G \big( (\alpha_{i})_{i \geq 0} \big)
$$

$(3)$ follows from \thref{classinteza1} and the fact that
$\mathcal{ZP}(A) = P(A) \simeq k^{(\mathbb{N})}$. For $(2)$ we use
\thref{classinteza1} again and the fact that $\Aut_{\rm Hopf}(A) =
\{u_{\alpha} \,\, | \,\, \alpha \in k^*\}$, where $u_{\alpha}$,
for $\alpha \in k^*$, is the automorphism given by $u_{\alpha}(Y)
= \alpha Y$. Thus, if $(\alpha_{i})_{i \geq 0}$, $(\beta_{i})_{i
\geq 0} \in k^{(\NN)}$ then $\sum_{i \geq 0} \alpha_{i} Y^{p^{i}}
\equiv \sum_{i \geq 0} \beta_{i} Y^{p^{i}}$ if and only if there
exists a pair $(\alpha, \beta) \in k^* \times k^*$ such that
$u_{\alpha} \big( \sum_{i \geq 0} \alpha_{i} Y^{p^{i}} \big) =
\beta^{2} \sum_{i \geq 0} \beta_{i} Y^{p^{i}}$. Since the equality
holds if and only if $\alpha^{p^{i}} \alpha_{i} = \beta^{2}
\beta_{i}$, for all $i \geq 0$, we obtain $(2)$.
\end{proof}

Using $(2)$ of \prref{exemplul2} we obtain:

\begin{corollary} \colabel{infinit1}
Let $k$ be a field of characteristic $p \geq 3$ and $(e_{j})_{j
\geq 0}$ the canonical $k$-basis of $k^{(\NN)}$, i.e. $e_{j} =
(\delta_{ij})_{i \geq 0}$, where $\delta_{ij}$ is Kronecker's
delta. Then $A_{(e_{j})}$, $j \geq 0$ defined by
\equref{relatiiex2} is an infinite family of non-isomorphic Hopf
algebras, i.e. $\textsc{C}\textsc{r}\textsc{p} (H_4, \, k[Y])$ is
an infinite set.
\end{corollary}

We now consider $A := k \lan y \, | \, y^p = 0 \ran$, the
$p$-dimensional Hopf algebra generated by a primitive element $y$.

\begin{proposition}\prlabel{exemplul3}
Let $k$ be a field of characteristic $p \geq 3$ and $A := k \lan y
\, ~|~ \, y^p = 0 \ran$ the $p$-dimensional Hopf algebra generated
by a primitive element $y$. Then:

$(1)$ Up to an isomorphism of Hopf algebras there exist exactly
two crossed products of Hopf algebras $A\#H_4$: $A \ot H_4$ and
$A_{4p}$, where $A_{4p}$ is the $4p$-dimensional Hopf algebra
generated by $g$ and $x$ subject to the relations:
\begin{equation} \eqlabel{relatiiex3}
x^{2p} = 0, \quad g^2 = 1, \quad x g = - g x
\end{equation}
and with the coalgebra structure given such that $g$ is a
group-like element and $x$ is $(1, g)$-primitive. Furthermore, we
have the following isomorphisms of groups:
$$
\Aut_{\rm Hopf} \big( A \otimes H_{4} \big) \simeq k^* \times k^*
\qquad {\rm and} \qquad \Aut_{\rm Hopf} \big( A_{4p} \big) \simeq
k^*
$$
$(2)$ There exists a bijection ${\mathcal H}^{2} (H_4, A) \cong
k$.
\end{proposition}

\begin{proof}
The proofs of $(2)$ follow exactly like the ones in
\prref{exemplul1}. For $(1)$ we only remark that, if $q \in k^*$
then $A_{(qy)}$ is the $4p$-dimensional quantum group generated by
$y$, $g$ and $x$, subject to the relations
$$
y^p = 0, \quad g^{2} = 1, \quad  x^{2} = qy, \quad  xg = - gx,
\quad gy = yg, \quad xy = yx
$$
which can be reduced to the ones in \equref{relatiiex3}.
\end{proof}

Finally, we consider $A := k \lan y \, | \, y^p = y \ran$, the
$p$-dimensional semi-simple Hopf algebra generated by a primitive
element $y$. For a positive integer $d$ we denote by $U_{d}(k)$
the set of $d$-th roots of unity in $k$.

\begin{proposition} \prlabel{exemplul4}
Let $k$ be a field of characteristic $p \geq 3$ and $A := k \lan y
\, | \, y^p = y \ran$ the $p$-dimensional Hopf algebra generated
by a primitive element $y$. Then:

$(1)$ A Hopf algebra $E$ is isomorphic to a crossed product $A \#
H_{4}$ if and only if $E \simeq A \ot H_{4}$ or $E \simeq
A_{(q)}$, for some $q \in k^*$, where $A_{(q)}$ is the
$4p$-dimensional quantum group generated by $g$ and $x$ subject to
the following relations
\begin{equation}\eqlabel{relatiiex4}
\quad g^2 = 1, \quad x^{2p} = q^{p-1}x^{2}, \quad xg = - gx
\end{equation}
and with the coalgebra structure given such that $g$ is a
group-like element and $x$ is $(1, g)$-primitive. Furthermore, we
have the following isomorphisms of groups:
$$
\Aut_{\rm Hopf} \big( A \otimes H_{4} \big) \simeq
\mathbb{F}_{p}^* \times k^* \qquad {\rm and} \qquad \Aut_{\rm
Hopf} \big( A_{(q)} \big) \simeq U_{2(p-1)}(k)
$$
for all $q \in k^*$.

$(2)$ There exists a bijection $\textsc{Crp} (H_4, \, A) \cong
k/_{\sim}$, where $\sim$ is the equivalence relation on $k$
defined by: $q \sim q'$ if and only if there exists $(\alpha,
\beta) \in \mathbb{F}_{p}^* \times k^*$ such that $\alpha q =
\beta^{2} q'$.

$(3)$ There exists a bijection ${\mathcal H}^{2} (H_4, A) \cong
k$.
\end{proposition}

\begin{proof}
For (1) we apply \coref{noucorH}. If $q \neq 0$ then $A_{(q y )}$
is generated by $y$, $g$ and $x$ subject to the relations
\begin{equation}
y^{p} = y, \quad g^2 = 1, \quad x^2 = qy, \quad xg = -gx, \quad gy
= yg, \quad xy = yx
\end{equation}
Since $y = q^{-1}x^{2}$, the independent relations are the ones in
\equref{relatiiex4}. For the automorphism groups we use
\reref{calculefec} and \leref{lemacumorfisme}. We have:
$$
\Aut_{\rm Hopf} \big( A \ot H_{4} \big) \simeq \Aut_{\rm Hopf}(A)
\times k^* \simeq \mathbb{F}_{p}^* \times k^*
$$
and
$$
\Aut_{\rm Hopf} \big( A_{(q)} \big) \simeq G(qy) = \{ (u, \beta)
\in \Aut_{\rm Hopf}(A) \times k^* \,\, | \,\, u(qy) = \beta^{2}qy
\} \simeq U_{2(p-1)}(k)
$$
where the isomorphism $G(y) \simeq U_{2(p-1)}(k)$ is given by
$G(qy) \ni (u, \beta) \mapsto \beta \in U_{2(p-1)}(k)$.

We now use \thref{classinteza1} for the classification part. From
\leref{lemacumorfisme} we have $\mathcal{ZP}(A) = P(A) = \{qy ~|~
q \in k\} \simeq k$, hence our claim in $(3)$. We look now at the
equivalence relations $\equiv$ on $\mathcal{ZP}(A)$. Let $q$, $q'
\in k$. Then $qy \equiv q'y$ if and only if there exists a pair
$(u, \beta) \in \Aut_{\rm Hopf}(A) \times k^*$ such that $u(qy) =
\beta^{2}q'y$. Following the description of $\Aut_{\rm Hopf}(A)$
in \leref{lemacumorfisme}, $qy \equiv q'y$ if and only if there
exists a pair $(\alpha, \beta) \in \mathbb{F}_{p}^* \times k^*$
such that $\alpha q = \beta^{2} q'$. This proves $(2)$.
\end{proof}

\prref{exemplul4} proves that the number of types of isomorphism
of all crossed products $ k \lan y \, | \, y^p = y \ran \# H_4$
depends heavily on the base field $k$. We illustrate this by the
following:

\begin{corollary} \colabel{nrtip1}
Let $p \geq 3$ be a prime number and $k = \mathbb{F}_p
(X_{1},X_{2}, \dots, X_{n}, \dots)$ the field of rational
functions in indeterminates $\{X_{i}\}_{i \geq 1}$ over the finite
field $\mathbb{F}_{p}$. Then $A_{(X_i)}$, $i \geq 1$, constructed
in \equref{relatiiex4} is an infinite family of non-isomorphic
$4p$-dimensional Hopf algebras.
\end{corollary}

\begin{proof}
Let $i$ and $j$ be two positive integers. Then $A_{(X_i)} \simeq
A_{(X_j)}$ as Hopf algebras if and only if $X_{i} \sim X_{j}$. We
claim that this is not the case if $i$ and $j$ are distinct.
Indeed, suppose $i \neq j$ and $X_{i} \sim X_{j}$. Then, by
\prref{exemplul4} (2), there exists a pair $(\alpha, \beta) \in
\mathbb{F}_{p}^* \times k^*$ such that $\alpha X_{i} = \beta^{2}
X_{j}$. Let $n \geq 1$ and $P$, $Q \in \mathbb{F}_{p}[X_{1},\dots,
X_{n}]$ such that $\beta = \frac{P}{Q}$. Then $\alpha X_{i} Q^{2}
= X_{j} P^{2}$. Considering both members of the equality as
polynomials in $X_{i}$ we obtain a contradiction, since the degree
of $\alpha X_{i} Q^{2}$ is odd, while the degree of $X_{j} P^{2}$
is even. Thus, $A_{(X_i)} \ncong A_{(X_j)}$, for any $ i \neq j
\geq 1$.
\end{proof}

\textbf{Acknowledgment.} We would like to thank to the referee for
his/her comments that substantially improved the first version of
this paper.


\begin{thebibliography}{99}

\bibitem{agorecia}
Agore, A.L., - Crossed product of Hopf algebras, {\sl Comm.
Algebra}, in press, arXiv:1203.2454.

\bibitem{ABM1}
Agore, A.L., Bontea, C.G. and Militaru, G. -- Classifying
bicrossed products of Hopf algebras, arXiv:1205.6110.

\bibitem{am1}
Agore, A.L. and Militaru, G. - Extending structures II: The
quantum version, {\sl J. Algebra}, {\bf 336} (2011), 321-341

\bibitem{AB}
Alperin, J.L. and Bell, R.R. - Groups and representations,
Springer-Verlag, New York, 1995.

\bibitem{AndCan}
Andruskiewitsch, N. -  Notes on Extensions of Hopf algebras, {\sl
Can. J. Math.} {\bf 48} (1996), 3--42.

\bibitem{AD}
Andruskiewitsch, N. and Devoto, J. - Extensions of Hopf algebras,
{\sl Algebra i Analiz} {\bf 7} (1995), 22--61.

\bibitem{AndN}
Andruskiewitsch, N. and Natale, S. - Examples of self-dual Hopf
algebras, {\sl J. Math. Sci. Univ. Tokyo}, {\bf 6}(1999), 181-215.

\bibitem{AMB10}
Ardizzoni, A., Beattie, M. and Menini, C. - Cocycle deformations
for Hopf algebras with a coalgebra projection, {\sl J. Algebra},
{\bf 324} (2010), 673--705.

\bibitem{AMSt}
Ardizzoni, A., Menini, C. and Stefan, D. - A monoidal approach to
splitting morphisms of bialgebras, {\sl Trans. AMS}, {\bf 359}
(2007), 991--1044

\bibitem{By1}
Byott, N. P.  - Cleft extensions of Hopf algebras, {\sl J.
Algebra}, {\bf 157}(1993), 405--429.

\bibitem{By2}
Byott, N. P.  - Cleft extensions of Hopf algebras II, {\sl Proc.
London Math. Soc.} {\bf 67}(1993), 277--304.

\bibitem{dnr}
Dascalescu, S., Nastasescu, C. and Raianu, S. - Hopf algebras. An
introduction, Marcel Dekker, New York, 2000.

\bibitem{DT}
Doi, Y. and Takeuchi, M. - Cleft comodule algebras for a
bialgebra, {\sl Comm. Algebra}, {\bf 14} (1986), 801 -- 818.

\bibitem{GV}
Garcia, G. A. and Vay, C. -  Hopf algebras of dimension 16, {\sl
Algebr. Represent. Theory},  {\bf 13} (2010), 383--405.

\bibitem{Ho}
Hofstetter, I - Extensions of Hopf algebras and thier
cohomological description, {\sl J. Algebra}, {\bf 164} (1994),
264--298.

\bibitem{MS}
Majid, S. and Soibelman, Ya. S. - Bicrossproduct structure of the
quantum Weyl group, {\sl J. Algebra} {\bf 163} (1994), 68 -- 87.

\bibitem{Ma1}
Masuoka, A. - Extensions of Hopf algebras, Trabajos de Matematica
41/99 Fa.M.A.F. (1999).

\bibitem{Ma2}
Masuoka, A. - Abelian and non-abelian second cohomologies of
quantized enveloping algebras, {\sl J. Algebra}, {\bf 320} (2008),
1--47.

\bibitem{Mo}
Molnar, R.K. - Semi-direct products of Hopf algebras, {\sl J.
Algebra} {\bf 47} (1977),  29 -- 51

\bibitem{M}
Montgomery, S. - Hopf algebras and their actions on rings, vol. 82
of CBMS Regional Conference Series in Mathematics, AMS,
Providence, Rhode Island (1993).

\bibitem{R}
Rotman, J.J. -  An introduction to the theory of groups. Fourth
edition. Graduate Texts in Mathematics 148, Springer-Verlag, New
York, 1995.

\bibitem{Sin}
Singer, W. - Extension theory for connected Hopf algebras, {\sl J.
Algebra}, {\bf 21}(1972), 1-16.

\bibitem{Sch2}
Schauenburg, P. - The structure of Hopf algebras with a weak
projection, {\sl Algebr. Represent. Theory} {\bf 3} (1999), 187
--211.

\bibitem{Sch93}
Schneider, H.-J. - Some remarks on exact sequences of quantum
groups, {\sl Comm. Algebra} {\bf 21}(1993), 3337--3358.

\bibitem{Sw68}
Sweedler, M.E. - Cohomology of algebras over Hopf algebras, {\sl
Trans. AMS}, {\bf 133} (1968), 205--239.

\end{thebibliography}
\end{document}